\documentclass[a4paper,leqno]{amsart}
\usepackage{amsfonts}
\theoremstyle{plain}
\begingroup

\newtheorem{theorem}{Theorem}[section]
\newtheorem{lemma}[theorem]{Lemma}
\newtheorem{proposition}[theorem]{Proposition}
\newtheorem{corollary}[theorem]{Corollary}
\newtheorem{remark}[theorem]{Remark}
\newtheorem{definition}[theorem]{Definition}
\endgroup

\theoremstyle{definition}
\theoremstyle{remark}

\mathchardef\emptyset="001F

\numberwithin{equation}{section}

\newcommand{\e}{\varepsilon}
\newcommand{\Om}{\Omega}
\newcommand{\Omk}{\Om\setmeno K}

\newcommand{\R}{{\mathbb R}}
\newcommand{\wto}{{\rightharpoonup}}
\newcommand{\setmeno}{\!\setminus\!}
\newcommand{\pinfty}{{+}\infty}

\newcommand{\huno}{{\mathcal H}^{1}}
\newcommand{\E}{{\mathcal E}}
\newcommand{\hb}{L^{1,2}}
\newcommand{\gdot}{{\dot{g}}}
\newcommand{\K}{{\mathcal K}(\overline\Om)}
\newcommand{\KA}{{\mathcal K}(\overline A)}

\title[Quasi-static growth of a brittle fracture]
{A  model for the quasi-static growth\\
of a brittle fracture:\\ existence
and approximation results
 }
\author[Gianni Dal Maso]{Gianni Dal Maso}
\address[Gianni Dal Maso]{SISSA, Via Beirut 2-4, 34014 Trieste, 
Italy}
\email[Gianni Dal Maso]{dalmaso@sissa.it}
\author[Rodica Toader]{Rodica Toader}
\address[Rodica Toader]{SISSA, Via Beirut 2-4, 34014 Trieste, 
Italy}
\email[Rodica Toader]{toader@sissa.it}
\begin{document}

\begin{abstract}
We give a precise mathematical formulation of a variational model for 
the irreversible quasi-static evolution of a brittle fracture proposed by 
G.A. Francfort and J.-J. Marigo, and based on Griffith's theory of 
crack growth.
In the two-dimensional case we prove an existence result for the 
quasi-static evolution and show that the total energy is an absolutely 
continuous function of time, although we can not exclude that
the bulk energy and the surface 
energy may present some jump discontinuities. This existence result 
is proved by a time discretization process, where at each step a 
global energy minimization is performed, with the constraint that the 
new crack contains all cracks formed at the previous time steps. This 
procedure provides an effective way to approximate the 
continuous time evolution. 
\end{abstract}
\maketitle
{\small

\bigskip
\keywords{\noindent {\bf Keywords:} variational models, 
energy minimization, free-discontinuity problems, 
crack propagation,
quasi-static evolution, brittle fracture, Griffith's criterion,
stress intensity factor.}


\bigskip
\subjclass{\noindent {\bf 2000 Mathematics Subject Classification:} 
35R35, 74R10, 49Q10, 35A35, 35B30, 35J25.}
}
\bigskip
\bigskip

\begin{section}{INTRODUCTION}

Since the pioneering work of A.~Griffith \cite{Gri}, the growth of a 
brittle fracture is considered to be the result of the competition 
between the energy spent to increase the crack and the corresponding
release of bulk energy. This idea is the basis of the celebrated 
Griffith's criterion for crack growth (see, e.g., \cite{SL}), and is used to 
study the crack propagation along a preassigned path. 
The actual path followed by 
the crack is often determined by using different criteria 
(see, e.g., \cite{ES}, \cite{SL}, \cite{SM}).

Recently G.A. Francfort and J.-J. Marigo \cite{FraMar3} proposed 
a variational model for the quasi-static growth of a brittle 
fracture, based on Griffith's theory, 
where the interplay between bulk and surface energy 
determines also the crack path.

The purpose of this paper is to give a precise mathematical 
formulation of a variant of this model in the {\it two-dimensional 
case\/}, 
and to prove an existence result for the {\it quasi-static evolution of a 
fracture\/} by using the {\it time discretization method\/} proposed in 
\cite{FraMar3}.

To simplify the mathematical description of the model, we consider 
only {\it linearly elastic homogeneous isotropic
materials\/}, with Lam\'e coefficients 
$\lambda$ and $\mu$. We restrict our analysis to
the case of an {\it anti-plane shear\/}, where 
the reference configuration is an infinite cylinder ${\Om{\times}\R}$, with 
$\Omega\subset \R^2$, and the displacement has the special form 
$(0,0,u(x_1,x_2))$ for every $(x_1,x_2,y)\in {\Om{\times}\R}$.
We assume also that the cracks  have the form
$K{\times}\R$, where $K$ is a compact set in $\overline\Om$. 
In this case the notions of bulk energy and surface energy 
refer to a finite portion of the cylinder determined by 
two cross sections separated by a unit distance.
The {\it bulk energy\/} is given by
\begin{equation}\label{bulk}
\frac{\mu}2\int_{\Om\setminus K} |\nabla u|^2dx\,,
\end{equation}
while the {\it surface energy\/} is given by
\begin{equation}\label{surface}
k\,\huno(K)\,,
\end{equation}
where $k$ 
is a constant which depends on the toughness of the material, and
$\huno$ is the {\it one-dimensional Hausdorff measure\/}, which 
coincides with the ordinary length in case $K$ is a rectifiable arc. For 
simplicity we take $\mu=2$ and $k=1$ in (\ref{bulk}) 
and (\ref{surface}).

We assume that $\Om$ is a {\it simply connected bounded open set\/} 
with a {\it Lipschitz boundary\/} $\partial\Om$, which, 
under these assumptions, is a simple closed curve. As in 
\cite{FraMar3}, we fix a subset $\partial_D\Om$ of $\partial\Om$, on 
which we want to  prescribe a {\it Dirichlet boundary condition\/} for $u$. We 
assume that $\partial_D\Om$ is composed of a {\it finite number of simple 
open arcs\/} with disjoint closures.

Given a function $g$ on $\partial_D\Om$, 
we consider the boundary condition 
$u=g$ on $\partial_D\Omk$. We can not prescribe a Dirichlet 
boundary condition on $\partial_D\Om\cap K$, because 
the boundary displacement is not 
transmitted through the crack, if the crack touches the boundary. Assuming 
that {\it the fracture is traction 
free\/} 
(and, in particular, without friction), 
the displacement $u$ in $\Omk$ is obtained by 
{\it minimizing (\ref{bulk}) under the boundary condition 
$u=g$ on $\partial_D\Omk$\/}. 
The {\it total energy} relative to the boundary displacement
$g$ and to the crack determined by $K$ is therefore
\begin{equation}\label{e0}
\E(g,K)=\min_v \Big\{\int_{\Om\setminus K}|\nabla 
v|^2dx+\huno(K) : v=g \hbox{ on }\partial_D\Om\setmeno K \Big\}\,.
\end{equation}
As $K$ is not assumed to be smooth,
we have to be careful in the precise mathematical formulation
of this minimum problem, which is given
at the beginning of  Section~\ref{conjugate}. The corresponding 
existence result is based on  some properties of the {\it Deny-Lions 
spaces\/}, that are described in Section~\ref{notation}.

In the theory developed in \cite{FraMar3} a crack with finite 
surface energy is any compact subset $K$ of $\overline\Om$ with 
$\huno(K)<\pinfty$. For technical reasons, that will be explained 
later, we propose a variant of this model, where we consider only 
cracks that are {\it compact and connected\/}
(these sets are called {\it continua\/}). There 
is no mechanical motivation for this choice, which is dictated by
mathematical convenience. Without this restriction some 
convergence arguments used in the proof of our existence result
are not justified
by the present development of 
the mathematical theories related to this subject.

We now describe our model of {\it quasi-static irreversible 
evolution of a fracture\/} under the action of a {\it time dependent 
boundary displacement\/} $g(t)$, $0\le t\le 1$. 
As usual, we assume that $g(t)$ can be extended to a function, still 
denoted by $g(t)$, which belongs to the Sobolev space $H^1(\Omega)$.
In addition, we assume that the function $t\mapsto g(t)$ is {\it 
absolutely continuous\/} from $[0,1]$ into $H^1(\Omega)$. Following 
the ideas of \cite{FraMar3}, given an initial (compact connected) crack $K_0$, 
with $\huno(K_0)<\pinfty$, we look for an {\it increasing family\/}
$K(t)$, $0\le t\le 1$, of (compact connected) cracks, 
with $\huno(K(t))<\pinfty$,
such that for any time
$t\in(0,1]$ the crack $K(t)$ {\it minimizes the total energy\/} 
$\E(g(t),K)$ among all (compact connected) cracks which contain 
all previous cracks $K(s)$, $s<t$. For $t=0$ we assume that
$K(0)$ minimizes $\E(g(0),K)$ among all (compact connected) cracks 
which contain $K_0$. 

This minimality condition for every time $t$ 
is inspired by Griffith's analysis of the energy balance. The constraint given
by the presence of the previous 
cracks reflects the {\it irreversibility of the evolution\/} 
and the {\it absence of a healing process\/}. In addition to this minimality 
condition we require also that
$\frac{d}{ds}\E(g(t),K(s))|_{s=t}=0$
for a.e.\ $t\in[0,1]$. In the special case $g(t)=t\,h$ for a given 
function $h\in H^1(\Om)$, we will see (Proposition~\ref{th}) that the last condition 
implies  the third condition  considered in Definition 2.9 of 
\cite{FraMar3}: $\E(g(t), K(t))\le \E(g(t), K(s))$ for $s<t$.

In Section~\ref{irrev} we prove the following existence result, where
$\K$ is the set of all continua $K$ contained 
in $\overline\Om$ and such that $\huno(K)<\pinfty$.
\eject
\begin{theorem}\label{kt0}
Let $g\in AC([0,1];H^1(\Om))$ and let $K_{0}\in\K$. Then 
there exists a function $K\colon [0,1]\to\K$ such that
\smallskip
\begin{itemize}
\item[(a)] 
\hfil $\displaystyle \vphantom{\frac{d}{ds}} 
K_0\subset K(s)\subset K(t)$ for $0\le s\le t\le 1$, \hfil
\item[(b)] \hfil $\displaystyle \vphantom{\frac{d}{ds}}
\E(g(0),K(0))\leq \E(g(0),K)
\quad\forall \,K\in\K,\,\ K\supset K_0$,\hfil
\item[(c)] \hfil $\displaystyle \vphantom{\frac{d}{ds}}
\hbox{for  }\, 0<t \le1\quad\E(g(t),K(t))\leq \E(g(t),K)
\quad\forall \, K\in\K,\,\ K\supset 
{\textstyle\bigcup_{s<t}K(s)}$,\hfil
\item[(d)]\hfil $\displaystyle \vphantom{\frac{d}{ds}}
t\mapsto \E(g(t),K(t)) \hbox{ is 
absolutely continuous on }[0,1]$, \hfil
\item[(e)]\hfil$\displaystyle
\frac{d}{ds}\E(g(t),K(s))\Big|\lower1.5ex\hbox{$\scriptstyle s=t$}=0
\quad \hbox{for a.e.\ }t\in[0,1]$.\hfil 
\end{itemize}
\smallskip
Moreover every function $K\colon[0,1]\to \K$ 
which satisfies (a)--(e) satisfies also
\begin{itemize}
\item[(f)] \hfil$\displaystyle\frac{d}{dt}\E(g(t),K(t))=
2\int_{\Om\setminus K(t)} \nabla u(t)\,\nabla \gdot(t)\,dx
\quad \hbox{ for a.e.\ } t\in[0,1]$,\hfil
\end{itemize}
where $u(t)$ is a solution of the minimum problem (\ref{e0}) which defines 
$\E(g(t),K(t))$, and $\gdot(t)$ is the time derivative of the function 
$g(t)$.
\end{theorem}

If $g(0)=0$ and $K_{0}\neq\emptyset$, conditions (a) and (b) imply 
that $K(0)=K_{0}$. If $g(0)=0$ and $K_{0}=\emptyset$, we can prove 
that there exists a solution of problem (a)--(e) with $K(0)=K_{0}$
(Remark~\ref{g0}).
We underline that, although we can not exclude that the surface energy 
$\huno(K(t))$ may present some jump discontinuities in time
(see \cite[Section 4.3]{FraMar3}), in our result
{\it the total energy is always an 
absolutely continuous function of time\/} by condition~(d). 

If $\partial_D\Om$ is sufficiently smooth, we can integrate  by parts 
the right hand side of (f) and, taking into account the Euler equation 
satisfied by $u(t)$, we obtain
\begin{equation}\label{energy}
\frac{d}{dt}\E(g(t),K(t))=
2\int_{\partial_D\Om\setminus K(t)} 
\frac{\partial u(t)}{\partial \nu}\, \gdot(t)\,d\huno
\quad \hbox{ for a.e.\ } t\in[0,1]\,,
\end{equation}
where $\nu$ is the outer unit normal to $\partial\Om$. Since the 
right hand side of (\ref{energy}) is the power of the force exerted 
on the boundary to obtain the displacement $g(t)$ on 
$\partial_D\Om\setmeno K(t)$,
equality (\ref{energy}) expresses the {\it conservation of energy\/} 
in our quasi-static model, where all kinetic effects are neglected.

The proof of this existence result is obtained by a time discretization process.
Given a time step $\delta>0$, for every integer $i\ge 0$ we set
$t_i^\delta:=i\delta$ and $g_i^\delta:=g(t_i^\delta)$. 
We define $K_i^\delta$, inductively, as a solution of the minimum problem
\begin{equation}\label{pidelta0}
\min_K\big\{ \E(g_i^\delta,K) : K \in \K,\
K \supset K_{i-1}^\delta \big\}\,,
\end{equation}
where we set $K_{-1}^\delta=K_0$.

Let $u_i^\delta$ be a solution of the minimum problem (\ref{e0}) 
which defines $\E(g_i^\delta,K_i^\delta)$. On 
$[0,1]$ we define  the step 
functions $K_\delta$ and $u_\delta$ by setting
$K_\delta(t):=K_{i}^\delta$ and
$u_\delta(t):=u_{i}^\delta$ for $t_i^\delta\le t<t_{i+1}^\delta$.

Using a standard monotonicity argument, we prove that there exists a 
sequence $(\delta_k)$ converging to $0$ such that, for every $t\in[0,1]$,
$K_{\delta}(t)$ converges to a continuum $K(t)$ in the Hausdorff 
metric as $\delta\to 0$ along this sequence. Then we can apply the results 
on the convergence of the solutions to mixed boundary 
value problems in cracked domains 
established in Section~\ref{convergence}, and we prove that, if $u(t)$ is a 
solution of the minimum problem (\ref{e0}) 
which defines $\E(g(t),K(t))$, then 
$\nabla u_{\delta}(t)$ converges to $\nabla u(t)$ strongly in 
$L^2(\Om,\R^2)$ as $\delta\to 0$ along the same sequence considered 
above.

The technical hypothesis that the sets $K_{\delta}(t)$ 
are connected plays a crucial role here. Indeed, if this hypothesis is 
dropped, the convergence in the Hausdorff 
metric of the cracks $K_{\delta}(t)$ to the crack $K(t)$
does not imply the convergence of the corresponding
solutions of the minimum problems, as shown by
many examples in homogenization theory that can be found, e.g.,
in \cite{Khr}, 
\cite{Mur},  \cite{Dam},  \cite{Att-Pic}, \cite{Cor}. These
papers show also that the hypothesis of connectedness would not 
be enough in dimension larger than two.

The results of Section~\ref{convergence} are similar to those 
obtained by A. Chambolle and F. Doveri in \cite{Ch-D} and by D. Bucur 
and N. Varchon in \cite{BucVar1} and \cite{BucVar}, 
which deal with the case of a pure Neumann 
boundary condition. Since we impose a Dirichlet boundary condition 
on $\partial_D\Om\setmeno K_{\delta}(t)$ and a Neumann 
boundary condition on the 
rest of the boundary, our results can not be deduced easily from 
these papers, so we give an independent proof, which uses the 
duality argument of \cite{BucVar}.

{}From this convergence result and from an approximation lemma with 
respect to the Hausdorff metric, proved in Section~\ref{Hausdorff}, we 
obtain properties (a), (b), (c), (e), and (f) in integrated form, 
which implies (d).

The time discretization process described above
turns out to be a useful tool for the proof of the 
existence of a solution $K(t)$ of the problem considered in 
Theorem~\ref{kt0}, and provides also an effective way for the numerical 
approximation of this solution (see \cite{BFM}), 
since many algorithms have been 
developed for the numerical solution of minimum problems of the form
(\ref{pidelta0}) (see, e.g., \cite{BZ}, \cite{Rich}, \cite{RM}, \cite{Bou}, 
\cite{Ch}, \cite{Bou-Cha}).

In Section \ref{tips} we study the motion of the tips of the time 
dependent
crack $K(t)$ obtained in Theorem~\ref{kt0}, assuming that, in some 
open interval $(t_0,t_1)\subset [0,1]$, the crack $K(t)$ has a fixed 
number of tips, that these tips move smoothly, and that their paths 
are simple, disjoint and do not intersect $K(t_0)$. We prove 
(Theorem~\ref{Griffith}) that in this case
{\it Griffith's criterion for crack growth\/} is valid in our model: the 
absolute value of the
{\it  stress intensity factor\/}
(see Theorem~\ref{Grisvard1} and Remark~\ref{stress}) 
of the solution $u(t)$ is less than or equal 
to $1$ at each tip for every $t\in (t_0,t_1)$, and it is equal to $1$ 
at a given tip for almost every instant $t\in (t_0,t_1)$ in which the tip moves with 
positive velocity.
\end{section}

\begin{section}{NOTATION AND PRELIMINARIES}\label{notation}

Given an open subset $A$ of $\R^2$, we say that $A$ has a Lipschitz 
boundary at a point $x\in\partial A$ if $\partial A$
is the graph of a Lipschitz function near $x$, in the sense that there exist
an orthogonal coordinate system
$(y_1,y_2)$, a rectangle ${\mathcal R}=(a,b)\times(c,d)$ containing $x$, 
and a Lipschitz function $\Phi\colon(a,b)\to(c,d)$, such that
$A\cap{\mathcal R}=\{y\in
{\mathcal R}: y_2<\Phi(y_1)\}$. The set of all these points $x$ 
is the {\it Lipschitz part of the boundary\/} and will be denoted by 
$\partial_{L}A$. 
If $\partial_{L}A=\partial A$, we say that $A$ has a Lipschitz boundary.

Besides the Sobolev space $H^1(A)$ we shall use also the 
{\it Deny-Lions space\/}
$\hb(A):=\{u\in L^2_{loc}(A)\;|\; \nabla u\in L^2(A;\R^2)\}$, which 
coincides with the space of all distributions $u$
on $A$ 
such that $\nabla u\in L^2(A;\R^2)$ (see, e.g., 
\cite[Theorem~1.1.2]{Ma}).
For the proof of the following result  we refer, e.g., to 
\cite[Section 1.1.13]{Ma}.
\begin{proposition}\label{closed}
The set $\{\nabla u: u\in \hb(A)\}$ is closed in $L^2(A;\R^2)$.
\end{proposition} 
Under some regularity assumptions on the boundary, the following
result holds.
\begin{proposition}\label{h1}
Let $u\in \hb(A)$ and $x\in\partial_L A$. 
Then there exists a neighbourhood $U$ of $x$ such that 
$u|_{A\cap U}\in H^1(A\cap U)$. In particular, if $A$ is bounded and 
has  a Lipschitz boundary, then $\hb(A)= H^1(A)$.
\end{proposition}
\begin{proof}
Let ${\mathcal R}$ be the rectangle given by the definition of Lipschitz 
boundary. It is easy to check that $A\cap{\mathcal R}$ has a Lipschitz 
boundary. The 
conclusion follows now from the Corollary
to Lemma 1.1.11 in \cite{Ma}.
\end{proof}


We recall some properties of the functions in the spaces 
$H^1(A)$ and $\hb(A)$, which are related to the notion of capacity. 
For more details we refer to \cite{EG}, \cite{HKM}, \cite{Ma}, and \cite{Zie}.

\begin{definition}
The capacity of an arbitrary subset $E$ of $\R^2$ is defined as
$$
{\rm cap} (E) :=\inf_{u\in{\mathcal U}_{E}}\left\{ 
\int_{\R^2}|\nabla u|^2dx+\int_{\R^{2}}|u|^{2}dx\right\},
$$
where ${\mathcal U}_{E}$ is the set of all functions $u\in H^1(\R^{2})$
such that $u\geq 1$ a.e.\ in a neighbourhood of~$E$.
\end{definition}
We say that a property is true {\it quasi-everywhere\/} on a set $E$, and 
write {\it q.e.\/}, if it holds on $E$ except on a set of capacity zero. As usual, 
the expression almost everywhere, abbreviated as a.e., refers to the 
Lebesgue measure.

A function $u\colon E\to\overline \R$ is said to be {\it 
quasi-continuous\/} on $E$ if 
for every 
$\e>0$ there exists an open set $U_{\e}$, with ${\rm cap}(U_{\e})<\e$, such 
that $u|_{E\setminus U_{\e}}$ is continuous on $E\setmeno U_{\e}$. 

It is known that every function $u\in\hb(A)$ has a {\it quasi-continuous 
representative\/} $\tilde u$, which is uniquely defined q.e.\ on 
$A\cup\partial_L A$, and satisfies
\begin{equation}\label{medie}
\lim_{\rho\to 0} \;
-\hskip-1.1em\int_{B_{\rho}(x)\cap A}|u(y)-\tilde u(x)|dy=0\qquad 
\hbox{for q.e.\ } x\in A\cup\partial_L A\,,
\end{equation}
where $-\hskip-.9em\int$ denotes the average and $B_\rho(x)$ is the 
open ball with centre $x$ and radius $\rho$.
If $u_n\to u$ strongly in $H^1(A)$, then a 
subsequence of $(\tilde u_n)$ converges to $\tilde u$ q.e. in 
$A\cup\partial_L A$. 
If $u,\, v\in \hb(A)$ and their traces coincide $\huno$-a.e.\ on 
$\partial_L A$ then $\tilde u$ and $\tilde v$ 
coincide q.e.\ on $\partial_L A$.

In the quoted books the quasi-continuous representatives are defined 
only on $A$. The straightforward definition of $\tilde u$ on $\partial_L A$ 
relies on the existence of extension operators for Lipschitz domains; 
the q.e.\ uniqueness of $\tilde u$ on $\partial_L A$ can be deduced from 
(\ref{medie}).
To simplify the notation we shall always identify each function $u\in 
\hb(A)$ with its quasi-continuous representative $\tilde u$.

Propositions~\ref{closed} and \ref{h1} imply the following result.
\begin{corollary}\label{complete} Assume that $A$ is connected, and 
let ${\Gamma}$ be a non-empty relatively open subset of $\partial A$ 
with ${\Gamma}\subset\partial_L A$. Then the space
$\hb_0(A,\Gamma):=\{u\in\hb(A):u=0\;\hbox{ q.e.\ on }{\Gamma}\}$ is a Hilbert space 
with the norm $\|\nabla u\|_{L^2(A;\R^2)}$. Moreover, if $(u_n)$ is a 
bounded sequence in $\hb_0(A,\Gamma)$, then there exist a subsequence, 
still denoted by $(u_n)$, and a function $u\in \hb_0(A,\Gamma)$ such 
that $\nabla u_n\wto \nabla u$ weakly in $L^2(A;\R^2)$.
\end{corollary}
\begin{proof} 
Let $(v_n)$ be  a Cauchy sequence in $\hb_0(A,\Gamma)$. 
We can construct an increasing sequence 
$(A_k)$ of connected open sets with 
Lipschitz boundary such that 
$A=\bigcup_k A_k$, and 
${\Gamma}=\bigcup_{k}(\partial A_k\cap\partial A)$.

By Proposition~\ref{h1} the functions $v_n$ belong to
$H^1(A_k)$ and $v_n=0$ q.e.\ on 
$\partial A_k\cap\partial A$. As 
$\huno(\partial A_k\cap\partial A)>0$ for $k$ large enough,
by the Poincar\'e inequality 
$(v_n)$ is a Cauchy sequence in $H^1(A_k)$, and therefore it converges 
 strongly in $H^1(A_k)$ to a function $v$ with $v=0$ q.e.\ on 
 $\partial A_k\cap\partial A$.  It is then easy to construct a function 
 $v\in \hb(A)$ such that $v=0$ q.e.\ on ${\Gamma}$ and 
$v_n\to v$ strongly in $H^1(A_k)$ for every $k$. As
 $(\nabla v_n)$ converges strongly in 
$L^2(A;\R^2)$, we conclude that 
$v_n\to v$ strongly in 
$\hb_0(A;\Gamma)$.

Let $(u_n)$ be a 
bounded sequence in $\hb_0(A,\Gamma)$. As in the previous part of the 
proof we deduce that $(u_n)$ is bounded in $H^1(A_k)$ for every $k$. 
By a diagonal argument we can prove that there exist a subsequence, 
still denoted by $(u_n)$, and a function $u\in \hb(A)$ such 
that $u_n\wto u$ weakly in $H^1(A_k)$ for every $k$. This implies 
$u=0$ q.e.\ on 
 $\partial A_k\cap\partial A$ for every $k$, hence 
 $u\in \hb_0(A,\Gamma)$. As
 $(\nabla u_n)$ is bounded in 
$L^2(A;\R^2)$, we conclude that 
$\nabla u_n\wto \nabla u$
weakly in $L^2(A;\R^2)$.
\end{proof}

Given a metric space $M$, the {\it Hausdorff distance\/} between two 
 closed subsets $K_{1},\, K_{2}$ of $M$ 
 is defined by
$$
d_{H}(K_{1},K_{2}):=\min\big\{1,\max\{ \sup_{x\in K_1}{\rm dist}(x,K_2),\sup_{y\in 
K_2}{\rm dist}(y,K_1)\}\big\}\,,
$$
with the usual conventions  ${\rm dist}(x,\emptyset)=\pinfty$ and 
$\sup\emptyset=0$, so that $d_{H}(\emptyset, K)=0$ if $K=\emptyset$ 
and $d_{H}(\emptyset, K)=1$ if $K\neq\emptyset$. 
We say that $K_n\to K$ in the Hausdorff metric if
${d_H(K_n,K)\to0}$. The following compactness theorem is 
well-known (see, e.g., 
\cite[Blaschke's 
Selection Theorem]{Rog}). 
\begin{theorem}\label{compactness}
Let $(K_n)$ be a sequence of closed subsets of a 
compact metric space $M$. Then there exists a subsequence which 
converges in the Hausdorff metric to a closed set 
$K\subset M$.
\end{theorem}

We recall that a {\it continuum\/} is a closed and connected set. The class 
of continua is stable under convergence in the Hausdorff metric. 
\begin{remark}\label{continua}{\rm 
Note that, if $K$ is a continuum in $\R^2$ and $x,\, y\in K$, then 
$\huno(K)\geq |x-y|$.
This implies that if $K,\, H$ are continua in $\R^2$, $K\subset H$, 
and $\huno(H\setmeno K)=0$, then either $H=K$, or $K=\emptyset$ and 
$H$ has only one element.}
\end{remark}

\end{section}

\begin{section}{PROPERTIES OF THE HARMONIC CONJUGATE}\label{conjugate}

Throughout the paper $\Om$ is a fixed bounded simply connected open 
subset of $\R^{2}$ with a Lipschitz boundary $\partial\Om$, which, 
under these assumptions, is a simple closed curve. We fix also a 
(possibly empty)
 subset $\partial_N\Om$ of $\partial\Om$, composed of a finite 
number of simple open arcs with disjoint closures, 
on which we shall impose a Neumann boundary 
condition. Let $\partial_S\Om$ be the finite set of the end-points of 
the arcs which compose $\partial_N\Om$, and let 
$\partial_D\Om:=\partial\Om\setmeno(\partial_N\Om\cup\partial_S\Om)$, 
which turns out to be the 
union of a finite number of simple open arcs. On this set
we want to impose a 
Dirichlet boundary condition.

Given a compact set $K$ in $\overline\Om$,
we consider the following boundary value problem:
\begin{equation}\label{*}
\left\{\begin{array}{ll}
\Delta u=0 & \hbox{in }\Omk\,,\\
\frac{\partial u}{\partial\nu}=0 & \hbox{on } 
\partial(\Omk)\cap(K\cup \partial_N\Om)\,.
\end{array}\right.
\end{equation}
By a solution of (\ref{*}) we mean a function $u$ which satisfies the 
following conditions:
\begin{equation}\label{**}
\left\{\begin{array}{l}
u\in \hb(\Omk)\,,\\
\displaystyle\int_{\Om\setminus K}\nabla u\,\nabla z\,dx=0\quad\forall z\in 
\hb(\Omk)\,,\ z=0\quad \hbox{q.e.\ on } \partial_D\Om\setmeno K\,.
\end{array}\right.
\end{equation}
Any solution $u$ of (\ref{**}) satisfies $\Delta u=0$ in the sense of 
distributions in $\Omk$ and, therefore, belongs to 
$C^\infty(\Omk)$.

Since no boundary condition is prescribed on $\partial_D\Om\setmeno K$, we do not 
expect a unique solution to problem (\ref{*}). Given $g\in \hb(\Omk)$,
we can prescribe the Dirichlet boundary condition 
\begin{equation}\label{**b}
u=g\quad\hbox{q.e.\ on }\partial_D\Om\setmeno K\,.
\end{equation}
It is clear that problem (\ref{**}) with the boundary condition
(\ref{**b}) can be solved separately in 
each connected component of $\Omk$. By Corollary~\ref{complete} and by 
the Lax-Milgram lemma there exists 
a unique solution in those components whose boundary meets 
$\partial_D\Om\setmeno K$, while on the other components the solution is 
given by an arbitrary constant. Thus the solution is not unique, if 
there is a connected component whose boundary does not meet 
$\partial_D\Om\setmeno K$. Note, however, that $\nabla u$ is always unique.
Moreover, the map $g\mapsto\nabla u$ is linear from $\hb(\Omk)$ 
into $L^2(\Omk;\R^2)$ and satisfies the estimate
$$
\int_{\Om\setminus K}|\nabla u|^2\,dx\le
\int_{\Om\setminus K}|\nabla g|^2\,dx\,.
$$

By standard arguments on the minimization of quadratic forms it is 
easy to see that $u$ is a solution of problem (\ref{**}) and 
satisfies the boundary condition (\ref{**b}) if and only if $u$ is a 
solution of the minimum problem 
\begin{equation}\label{minpb}
\min_{v\in {\mathcal V}(g,K)}\int_{\Om\setminus K}|\nabla v|^2\,dx\,,
\end{equation}
where 
\begin{equation}\label{vgk}
{\mathcal V}(g,K):=\{v\in\hb(\Om\setmeno K): v=g\quad \hbox{q.e.\ on }
\partial_D\Om\setmeno K\}\,.
\end{equation}
Throughout the paper, given a function $u\in \hb(\Omk)$, we always 
extend $\nabla u$ to $\Om$ by setting $\nabla u=0$ a.e.\ on $K$. Note that, 
however, 
 $\nabla u$ is the  distributional gradient  of $u$ only in $\Omk$, and, in 
 general, it does not coincide in $\Om$ with the gradient 
 of an extension of $u$.

In the next lemma we prove, under some additional regularity 
assumptions, that every solution $u$ of (\ref{*}) 
has a harmonic conjugate $v$ on $\Omk$ which vanishes on $K$ and is 
constant on each connected component of $\partial_N\Om$. Let $R$ be the 
rotation on $\R^2$ given by $R(y_{1},y_{2})=(-y_{2},y_{1})$.

\begin{lemma}\label{unoreg} Assume that $\partial_N\Om$ is a 
manifold  of 
class $C^\infty$. Let $K$ be a continuum that can be written in the 
form $K:=\overline A\cap\overline\Om$,
where $A$ is an open set with 
a $C^\infty$ boundary, $\partial A\cap\partial_S\Om=\emptyset$, 
$\partial A\cap\partial\Om$ has a finite 
number of points,  and $\partial A$ meets $\partial_N\Om$ forming angles 
different from $0$ and $\pi$.
Let $u$ be a solution of problem (\ref{**}), let 
$E:=\overline\Om\setmeno \overline{\partial_D\Om\setmeno K}$, and let
$F:=\partial A\cap\partial_D\Om$. Then
there exists a function $v\in H^1(\Om)\cap C^0(E)\cap C^\infty(\Omk)$, 
such that 
$\Delta v=0$ in $\Omk$, $\nabla v=R\,\nabla u$ in $\Omk$, $v=0$  
on $K\setmeno F$, 
and $v$ is constant on each connected component of $\partial_N\Om$.
\end{lemma}
\begin{proof} 
If $C$ is a connected component of $\Omk$ whose boundary does not 
meet $\partial\Omk$, then $C$ is simply connected. 
Indeed, if $U$ is a bounded open set with 
$\partial U\subset C$, then $ U\subset C$,
since otherwise $K\cap  U\neq\emptyset$ and 
$K\setmeno\overline U\neq\emptyset$, which contradicts the 
fact that $K$ is connected.
Therefore in this case $K\cap\partial C=\partial C$ is connected.

Let $C$ be a connected component of $\Omk$ whose boundary
meets $\partial\Omk$.
If $K\cap\partial \Om\neq\emptyset$, by applying the previous result to 
the continuum $K':=K\cup\partial\Om$, we obtain that $C$ is 
simply connected. By our regularity assumptions $\partial C$ is a 
closed simple curve. We want to prove that $K\cap\partial C$ is connected.
If not, there exist $x,\, y\in\partial\Om\cap\partial C$ such that
$K$ intersects both connected components of $\partial C\setmeno\{x,y\}$.
Let $\Gamma$ be a simple arc connecting $x$ and $y$, whose interior 
points lie in $C$. As $\Om$ is simply connected, $\Gamma$ divides 
$\overline\Om$ into two relatively open connected components, each of which 
contains one of the connected components of $\partial 
C\setmeno\{x,y\}$, and hence intersects $K$. As $K$ does not meet 
$\Gamma$, this contradicts the fact that $K$ is connected.
Therefore $K\cap\partial C$ is connected.

If $K\cap\partial\Om=\emptyset$, then $\Omk$ has exactly one 
connected component $C$ whose boundary 
intersects $\partial\Om$. Let $H$ be the connected component of 
$\partial K$ which contains $K\cap\partial C$. By our regularity 
assumption, it is easy to see that $K\cap\partial C$ is closed and 
open in $H$, so that $K\cap\partial C$ coincides with $H$, and hence 
it is connected.


Let us prove that in each connected component $C$ of $\Omk$ there exists 
a harmonic function $v$ such that $\nabla v=R\,\nabla u$. As the
differential form $-D_2u\,dx_1+D_1u\,dx_2$ is closed, the existence of 
$v$ is trivial if $C$ is simply connected. The only non-trivial case 
is when $K\cap\partial\Om=\emptyset$ and $C$ is the unique 
connected component of $\Omk$ whose boundary 
intersects $\partial\Om$.
Let $ U$ be an open set with smooth boundary such that 
$K\subset U\subset\subset\Om$. Since $u$ is a solution of 
(\ref{**}), integrating by parts we obtain that 
$\int_{\partial  U}\frac{\partial u}{\partial\nu}d\huno=0$.
This implies 
that the integral on $\partial U$ of the differential form 
$-D_2u\,dx_1+D_1u\,dx_2$ vanishes, and proves the existence of the 
harmonic conjugate also in this case.

Let $\tau$ and $\nu$ be  smooth tangent and normal vectors to 
$\partial K\cap\Om=\partial A\cap\Om$.  
By classical regularity results for solutions to Neumann problems 
on domains with corners (see \cite[Theorem 4.4.3.7]{Gri1}), 
we have that  $\nabla u$ has a continuous extension to 
$\overline{\Omk}\setminus\overline{\partial_D\Om\setmeno K}$
and $\frac{\partial u}{\partial\nu}=0$ on $\Om\cap\partial A$ and on
$\partial_N\Om\setmeno K$. 
This implies that $\nabla v$ and $v$ have a continuous extension to 
$\overline{\Omk}\setminus\overline{\partial_D\Om\setmeno K}$ and
$\frac{\partial v}{\partial\tau}=0$ on $\Om\cap\partial A$ and on
$\partial_N\Om\setmeno K$, hence $v$ is constant on 
$(K\cap\partial C) \setmeno F$ 
for every connected component $C$ of $\Omk$ (as 
$K\cap\partial C$ is connected, so is 
$(K\cap\partial C) \setmeno F$, which is 
obtained by possibly removing some of its end-points). 
We can choose the integration constant in each connected component 
so that $v=0$ on $K\cap\partial C$, 
therefore if we define $v=0$ in $\overline\Om\cap A$ we obtain that 
$v$ belongs to $H^1(\Om)$, is continuous in $E$, and $v=0$ on $K\setmeno F$.
Moreover $v$ is harmonic in $\Omk$ and $\nabla v=R\,\nabla u$ in $\Omk$.
As $\frac{\partial v}{\partial\tau}=0$ on
$\partial_N\Om\setmeno K$, $v$ is constant on each connected component of 
$\partial_N\Om\setmeno K$. As $v=0$ on $K\setmeno F$, the continuity of 
$v$ implies that $v$ is constant on each connected component of 
$\partial_N\Om$.
\end{proof}

To extend the previous result to the case of a general continuum $K$, we
use the following lemma, which is established in its more general form 
in view of the applications in the next sections.

\begin{lemma}\label{glinf} Suppose that $\Om$ is the union of an 
increasing 
sequence $(\Om_n)$ of open sets.
Let $(K_n)$ be a sequence of compact sets contained in 
$\overline\Om$ which converges to a compact set $K$ in the Hausdorff 
metric.
Let $u_n\in \hb (\Om_n\setmeno K_n)$ be a sequence such that 
$\|\nabla u_n\|_{L^2(\Om_n\setminus K_n;\R^2)}$ is bounded. Then there 
exist a subsequence, still denoted by $(u_n)$, and a function 
$u\in \hb (\Omk)$, such that $\nabla u_n\wto \nabla u$ 
weakly in $L^2(U;\R^2)$ for every open set $U\subset\subset\Omk$.
If, in addition, ${\rm meas}(K_n)\to{\rm meas}(K)$, then $\nabla u_n\wto 
\nabla u$ weakly in $L^2(\Om;\R^2)$, where we set $\nabla u_n=0$ on 
$(\Om\setmeno\Om_n)\cup K_n$ and $\nabla u=0$ on $K$. 
\end{lemma}
\begin{proof}
For every open set $U\subset\subset \Omk$ we have 
$U\subset\subset \Om_n\setmeno K_n$ for $n$ large enough.
Since the sequence $(\nabla u_n)$ is bounded in $L^2(U;\R^2)$, 
there exists a
subsequence, still denoted by $(u_n)$, such that $(\nabla u_n)$
converges weakly in $L^2(U;\R^2)$ to some function 
$\varphi\in L^2(U;\R^2)$.
 Since the space
$\{\nabla v: v\in \hb (U)\}$ is closed in $L^2(U;\R^2)$, we 
conclude that there exists $u\in \hb (U)$ such that 
$\nabla u=\varphi$ a.e.\ in $U$. Repeating the same argument for 
all $U\subset\subset \Omk$, we construct $u\in \hb (\Omk)$
such that $\nabla u=\varphi$ a.e.\ in $\Omk$.

Assume now that ${\rm meas}(K_n)\to{\rm meas}(K)$ and let $\psi\in 
L^2(\Om;\R^2)$. For every $\e>0$ there exists $\delta>0$ such that 
$\int_A|\psi|^2\,dx<\e^2$ for ${\rm meas}(A)<\delta$. Let 
$U\subset\subset \Omk$ be an open set such that 
${\rm meas}((\Om\setmeno K)\setmeno U))<\delta$. As 
$U\subset\subset \Om_n\setmeno K_n$ for $n$ large enough, we have also 
${\rm meas}((\Om_n\setmeno K_n)\setmeno U))<\delta$. Then 
$$
\Big|\int_{\Om}(\nabla u_n-\nabla u)\cdot\psi\,dx\Big|\leq
\Big|\int_{U}(\nabla u_n-\nabla u)\cdot\psi\,dx\Big| + c_1\e+c_2\e\,,
$$
where $c_1$ is an upper bound for 
$\|\nabla u_n\|_{L^2(\Om_n\setminus K_n;\R^2)}$ and 
$c_2:=\|\nabla u\|_{L^2(\Om\setminus K;\R^2)}$.
{}From the previous part of the lemma 
$\limsup_{n}|\int_\Om(\nabla u_n-\nabla u)\cdot\psi\,dx|\leq 
c_1\e+c_2\e$ and the conclusion follows from the arbitrariness of~$\e$.
\end{proof}

We are now in a position to prove the main result of this section.

\begin{theorem}\label{uno} Let $K$ be a continuum in $\overline\Om$ 
and let $u$ be a solution of problem (\ref{**}). Then
there exists a function $v\in H^1(\Om)\cap C^{\infty}(\Omk)$ such that 
$\Delta v=0$ in $\Omk$,
$\nabla v=R\,\nabla u$ in $\Omk$,
$v=0$ q.e.\ on $K$, and $v$ is constant q.e.\ on each connected 
component of $\partial_N\Om$.
\end{theorem}
\begin{proof} We can express $\Om$ as the union of an increasing sequence 
$(\Om_j)$ of simply connected open sets, with 
$\partial_D\Om\cup\partial_S\Om\subset\partial\Om_j$, such that
$\partial_N\Om_j:=\partial\Om_j\setmeno
(\partial_D\Om\cup\partial_S\Om)$
is a finite union of regular open arcs of class $C^\infty$ that are 
in one-to-one correspondence with the connected components of
$\partial_N\Om$, so that corresponding arcs have the same end-points.

We can also write $K$ as the intersection of a 
decreasing sequence $(A_j)$ of 
open sets with  boundary of class $C^\infty$ such that 
$K_j:=\overline A_j\cap\overline\Om_j$ is connected, 
$\partial A_j\cap\partial_S\Om=\emptyset$, 
$\partial A_j\cap\partial\Om_j$ has a finite number of points,
and $\partial A_j$ meets $\partial_N\Om_j$ forming
angles different from $0$ and $\pi$.

As $u\in\hb(\Om_j\setmeno K_j)$, there exists a solution $u_j$ to the problem 
\begin{equation}\label{**j}
\left\{\begin{array}{l}
u_j\in \hb(\Om_j\setmeno K_j)\,,\quad u_j=u\quad\hbox{q.e.\ on }
\partial_D\Om\setmeno K_j,\\
\displaystyle\int_{\Om_j\setminus K_j}\nabla u_j\,\nabla z\,dx=0\quad\forall z\in 
\hb(\Om_j\setmeno K_j)\,,\ z=0\quad\hbox{q.e.\ on }
\partial_D\Om\setmeno K_j\,.
\end{array}\right.
\end{equation}
Using $u_j-u$ as test function in (\ref{**j}), we obtain that the 
norms $\|\nabla u_j\|_{L^2(\Om_j\setminus K_j)}$ are uniformly bounded. By 
Lemma~\ref{glinf}, there exists $u^*\in\hb(\Omk)$ such that,  up to a 
subsequence,  $(\nabla u_j)$ converges 
to $\nabla u^*$ weakly  in $L^2(\Om;\R^2)$.

Let us prove that 
\begin{equation}\label{uu*}
\nabla u^*=\nabla u\quad\hbox{ a.e.\ in }\Omk\,.
\end{equation}
To this end it is enough to construct a solution $w$ of (\ref{**}) such that 
$\nabla u^*=\nabla w$ a.e.\ in $\Omk$ and 
$w=u$ q.e.\ on $\partial_D\Om\setmeno K$. 

Let $C$ be a connected component of $\Omk$ whose boundary 
meets $\partial_D\Om\setmeno K$, and let $x\in C$.  Given $\e>0$ 
small enough, let 
$N^{\e}:=\{x\in\R^2:{\rm dist}(x,\partial_N\Om\cup K)\leq\e\}$ and 
let $C^{\e}$ be the connected component of $C\setmeno N^{\e}$ 
containing $x$. 

Let ${\Gamma}^{\e}$ be the relative interior of 
$\partial C^{\e}\cap \partial_D\Om$ in $\partial C^{\e}$. As 
$\partial C$ meets $\partial_D\Om\setmeno K$, for $\e$ small enough 
${\Gamma}^{\e}\neq\emptyset$.  Using Corollary~\ref{complete} we 
deduce that there exists $w\in \hb(C^{\e})$ such that $w=u$ q.e.\ 
on $ {\Gamma}^{\e}$ and, up to a subsequence, $\nabla u_j\wto\nabla 
w$ weakly in $L^2(C^{\e};\R^2)$. Hence $\nabla w=\nabla u^*$ a.e.\ 
in $C^{\e}$. Since $\e>0$ is arbitrary, we can construct a function 
$w\in\hb(C)$ such that $w=u$ q.e.\ on $\partial 
C\cap(\partial_D\Om\setmeno K)$ and $\nabla w=\nabla u^*$ a.e.\ in $C$. 

Let now $C$ be a connected component of $\Omk$ whose boundary does not 
meet  $\partial_D\Om\setmeno K$ and let $x\in C$. It is easy to see 
that $C$ is the union of the increasing sequence $(C_j)$ of the 
connected components of $\Om_j\setmeno K_j$ containing $x$. Since 
$\partial C_j$ does not meet $\partial_D\Om\setmeno K_j$ and $u_j$ is 
a solution of (\ref{**j}), we conclude that $\nabla u_j=0$ a.e.\ in 
$C_j$, hence $\nabla u^*=0$ a.e.\ in $C$. 
We define $w:=0$ in these connected components.

In this way we have defined a function $w$ on each connected 
component of $\Omk$ such that $w\in \hb(\Omk)$, $w=u$ q.e.\ on 
$\partial_D\Om\setmeno K$,  $\nabla w=\nabla u^*$ a.e.\ in $\Omk$, 
and hence $\nabla u_j$ converges to $\nabla w$ weakly in 
$L^2(\Om;\R^2)$.

To conclude the proof of (\ref{uu*}) it is 
enough to show that $w$ is a solution of (\ref{**}).
Let $z\in\hb(\Omk)$ with $z=0$ q.e.\ on $\partial_D\Omk$. As $z\in 
\hb(\Omk_j)$ and $z=0$ q.e.\ on $\partial_D\Om\setmeno K_j$, 
we can use $z$ as test function in (\ref{**j}). Then passing to the limit as 
$j\to\infty$ we obtain (\ref{**}), and the proof of (\ref{uu*}) is 
complete.

Let $v_j\in H^1(\Om_j)$, with $v_j=0$ 
q.e.\ on $K_j$,  be the harmonic conjugate of $u_j$ on 
$\Om_j\setmeno K_j$ given by Lemma~\ref{unoreg}. 
As $\nabla v_j=R\,\nabla u_j$ a.e.\ on $\Om_{j}\setmeno K_j$ and 
$\nabla v_j=0$ a.e.\ on $K_j$, we deduce that $(\nabla v_j)$ is bounded in 
$L^2(\Om;\R^2)$. Since $v_j$ is constant on each connected component 
of $\partial_N\Om_j$, our hypotheses on the sets $\Om_j$ imply that 
this function can be extended to a function, still denoted by $v_j$, which 
belongs to $H^1(\Om)$ and is 
locally constant on $\overline\Om\setmeno\overline\Om_j$. In 
particular we have that $v_j$ is constant on each connected component 
of $\partial_N\Om$ and $(\nabla v_j)$ is bounded in $L^2(\Om;\R^2)$.

Assume that $K$ has more than one point. Then 
$\liminf_{j}{\rm diam}(K_j)>0$; 
since the sets $K_j$ are connected, we obtain also
$\liminf_{j}{\rm cap}(K_j)=:\delta>0$. As $v_j=0$ q.e.\ on $K_j$,
using the Poincar\'e inequality (see, e.g. \cite[Corollary 4.5.3]{Zie})
it follows that $(v_j)$ is bounded in $H^1(\Om)$, hence there
exists a function $v\in H^1(\Om)$ such that, up to a subsequence,  
$v_j\wto v$ weakly in $H^1(\Om)$. 
This implies that $\nabla v=R\,\nabla u$ 
a.e.\ in $\Omk$, that $v=0$ q.e.\ in 
$K$, and that $v$ is constant q.e.\ on each connected component of 
$\partial_N\Om$.

Assume now that $K$ has only one point, and let $v_j^\Om$ be the 
mean value of $v_j$ in $\Om$. By the Poincar\'e inequality, the 
sequence $(v_j-v_j^\Om)$ is bounded in $H^1(\Om)$, hence a 
subsequence converges to a function $v$ weakly in $H^1(\Om)$. It is 
clear that  $\nabla v=R\,\nabla u$ a.e.\ in $\Omk$ and that 
$v$ is constant q.e.\ on 
each connected component of 
$\partial_N\Om$. The condition $v=0$ q.e.\ on $K$ is trivial in this 
case.
\end{proof} 

\begin{theorem}\label{due}
Let $K$ be a continuum in $\overline\Om$ and let $v\in 
H^1(\Om)\cap C^\infty(\Omk)$ be a function such that $\Delta v=0$ in 
$\Omk$, $v=0$ q.e.\ on $K$, and $v$ is constant q.e.\ on each connected 
component of $\partial_N\Om$. Assume also that the differential form 
$-D_2v\,dx_1+D_1v\,dx_2$ is exact in $\Omk$. Then there exists 
a solution $u$ of 
(\ref{**}) such that $R\,\nabla u=\nabla v$ a.e.\ in $\Omk$.
\end{theorem}
Note that the assumption on the differential form is always satisfied  if 
$K\cap\partial\Om\neq\emptyset$, since, 
in this case, each connected component of 
$\Omk$ is simply connected (see the proof of Lemma~\ref{unoreg}). 
\begin{proof}
Since the differential form 
$-D_2v\,dx_1+D_1v\,dx_2$ is exact in $\Omk$, it follows that there 
exists a function $u\in C^{\infty}(\Omk)$ such that $R\,\nabla u=\nabla v$ 
 in $\Omk$ and 
$\Delta u=0$ in $\Omk$. 
As $\nabla v\in L^2(\Omk;\R^2)$, we have also 
$\nabla u\in L^2(\Omk;\R^2)$, hence $u\in\hb(\Omk)$.
It remains to prove that $u$ is a solution of (\ref{**}). 

Let us prove that, if $K$ has more than one point and $K$ meets the 
closure $\Gamma$ of a connected component of $\partial_N\Om$,
then $v=0$ q.e.\ on $\Gamma$. If this is not true, there exists a 
constant $c\neq 0$ such that $v=c$ q.e.\ on $\Gamma$. Let us fix 
$x\in K\cap\Gamma$. For almost every sufficiently small $\rho>0$ we 
have $\Gamma\cap\partial B_\rho(x)\neq \emptyset$ and 
$K\cap\partial B_\rho(x)\neq \emptyset$, hence $v$ takes the values 
$c$ and $0$ in two distinct points of $\partial B_\rho(x)$. This 
implies
$$
\int_{\partial B_\rho(x)}|\nabla v|^2\,d\huno\ge 
\frac{c^2}{2\pi\rho}\,,
$$
which yields $\int_\Om|\nabla v|^2\,dx=+\infty$, in contradiction with 
our hypothesis.

Therefore, if $K$ has more than one point, then $v=0$ q.e.\ on all 
connected components of $\partial_N\Om$ whose closure meets $K$. Since 
$v$ is constant q.e.\ on the other connected components
of $\partial_N\Om$, we can apply
\cite[Theorem~4.5]{HKM} to a suitable extension of $v$ and we can 
construct a sequence of functions $v_n\in C^{\infty}(\R^2)$, converging 
to $v$ in $H^1(\Om)$, such that each $v_n$ vanishes  in a 
neighbourhood of $K$ and is constant in a neighbourhood of each 
connected component of $\partial_N\Om$. 

Let $z\in \hb(\Omk)$ 
with $z=0$ q.e.\ on $\partial_D\Om\setmeno K$. 
As ${\rm div}(R\,\nabla v_n)=0$ in $\R^2$ and $R\,\nabla v_n=0$ 
near $K\cup\partial_N\Om$, we 
have
$$
\int_{\Om\setminus K}R\,\nabla v_n\,\nabla z\,dx=
\int_{\partial_D\Om\setminus K} R\,\nabla v_n z\,\nu\,d\huno =0\,.
$$ 
Passing to the limit as $n\to\infty$,
we obtain
$$
\int_{\Om\setminus K}\nabla u\,\nabla z\,dx=
-\int_{\Om\setminus K}R\,\nabla v\,\nabla z\,dx=0\,,
$$
showing thus that $u$ is a solution of (\ref{**}).
\end{proof}
\end{section}

\begin{section}{CONVERGENCE OF MINIMIZERS}\label{convergence}

In this section we prove the convergence of the minimum points of problems 
(\ref{minpb}) corresponding to a sequence $(K_n)$ of continua which 
converges  
in the Hausdorff metric to a 
continuum $K$ such that ${\rm meas}(K_n)\to{\rm meas}(K)$.

 \begin{theorem}\label{convsol}
Let $(g_n)$ be a sequence in $H^1(\Om)$ which converges to a 
function $g$ strongly in $H^1(\Om)$, and let $(K_n)$ be a sequence of  
continua in $\overline\Om$ 
which converges to a continuum $K$ in the Hausdorff metric.  Assume 
that ${\rm meas}(K_n)\to{\rm meas}(K)$. 
Let 
$u_n$ be a solution 
of the minimum problem 
\begin{equation}\label{Pn}
\min_{v\in{\mathcal V}(g_n,K_n)}
\int_{\Om\setminus K_{n}}|\nabla v|^2\,dx\,,
\end{equation}
and let $u$ be  a solution of the minimum problem 
\begin{equation}\label{P}
\min_{v\in{\mathcal V}(g,K)} \int_{\Om\setminus K}|\nabla v|^2\,dx\,.
\end{equation}
Then $\nabla u_n\to \nabla u$ strongly in $L^2(\Om;\R^2)$.
\end{theorem}
\begin{proof}
Note that $u$ is a minimum point  of (\ref{P}) if and only if $u$ satisfies 
(\ref{**}) and (\ref{**b}); analogously, $u_n$ is a minimum point  of 
(\ref{Pn}) 
if and only if  $u_n$ satisfies (\ref{**}) and (\ref{**b}) with $K$ 
and $g$ replaced 
by $K_n$ and $g_n$.

Taking $u_n-g_n$ as test function in the equation satisfied by $u_n$, we
prove that the sequence
$(\nabla u_n)$ is bounded in $L^2(\Om;\R^2)$. 
By Lemma~\ref{glinf},
there exists a function $u^*\in \hb(\Omk)$
such that, passing to a subsequence, $\nabla u_n\wto \nabla u^*$ 
weakly in $L^2(\Om;\R^2)$.

Let $v_n$ be the harmonic conjugate of $u_n$ given by Theorem~\ref{uno}. 
Then
$(\nabla v_n)$ is bounded in $L^2(\Om;\R^2)$, $v_n=0$ q.e.\ on $K_n$, 
and $v_n$ is constant q.e.\ on each connected component of 
$\partial_N\Om$.

If $K$ has only one point, using Propositions~\ref{closed} and 
\ref{h1}, we  
obtain that $\nabla v_n\wto \nabla v$ weakly in $L^2(\Om;\R^2)$ for 
some $v\in H^1(\Om)$, and by the Poincar\'e inequality there are 
constants $c_n$ such that $v_n-c_n\wto v$ weakly in $H^1(\Om)$.  This 
implies  that  $v$ is 
constant q.e.\ on each connected component of $\partial_N\Om$. 
Moreover $v=0$ q.e.\ 
on $K$ since ${\rm cap}(K)=0$.

If $K$ has more than one point, we have $\liminf_{n}{\rm diam}(K_n)>0$; 
since 
the sets $K_n$ are connected, we obtain also
$\liminf_{n}{\rm cap}(K_n)=:\delta>0$. As $v_n=0$ q.e.\ on $K_n$,
using the Poincar\'e inequality (see, e.g., \cite[Corollary 4.5.3]{Zie})
it follows that $(v_n)$ is bounded in $H^1(\Om)$. Hence 
there exists a function $v\in H^1(\Om)$ such that 
$v_n\wto v$ weakly in $H^1(\Om)$.  This implies  that $v$ is constant on 
each connected component of $\partial_N\Om$.

Let us fix an open ball $B$ containing $\overline\Om$. Using the same 
extension operator we can construct extensions of $v_n$ and $v$, 
still denoted by  $v_n$ and $v$, such that $v_n,\,v\in H^1_0(B)$, 
$v_n\wto v$ weakly in $H^1(B)$. 

Given an open set $A\subset B$, any function $z\in H^1_0(A)$ will 
be extended to a function $z\in H^1_0(B)$ by setting $z:=0$ q.e.\ in 
$\overline B\setmeno A$. By \cite[Theorem 4.5]{HKM} we have
\begin{equation}\label{h10}
H^1_0(A)=\{z\in H^1(B):z=0\quad\hbox{q.e.\ on }
\overline B\setmeno A\}\,.
\end{equation}
Since the complement of $B\setmeno K_n$ has two connected 
components, 
from the results of \cite{Sve} and \cite{Buc} we deduce that, for 
every $f\in L^2(B)$, the solutions $z_n$ of the Dirichlet problems 
$$
z_n\in H^1_0(B\setmeno K_n)\qquad \Delta z_n=f\quad\hbox{in }
B\setmeno K_n
$$
converge strongly in $H^1_0(B)$ to the solution $z$ of the Dirichlet 
problem 
$$
z\in H^1_0(B\setmeno K)\qquad \Delta z=f\quad\hbox{in }
B\setmeno K\,.
$$
This implies (see, e.g., \cite[Theorem~3.33]{Att})
that, in the space $H^1_0(B)$, the subspaces 
$H^1_0(B\setmeno K_n)$ converge to  the subspace 
$H^1_0(B\setmeno K)$ 
in the sense of Mosco (see \cite[Definition 1.1]{Mos}). 

 Since $v_n\in H^1_0(B\setmeno K_n)$ by (\ref{h10}),
and  $v_n\wto v$ weakly in $H^1(B)$, from the 
convergence in the sense of Mosco we deduce that 
$v\in H^1_0(B\setmeno K)$, hence $v=0$ q.e.\ on $K$ by (\ref{h10}).

As $v_n$ is harmonic in $\Omk_n$ and 
$\nabla v_n=R\,\nabla u_n$ in 
$\Omk_n$, we deduce that $v$ is harmonic in 
$\Omk$ and $\nabla v=R\,\nabla u^*$  in 
$\Omk$. Let us consider a smooth 
closed curve $\Gamma$ in $\Omk$. Since, for $n$ large enough, the 
functions $v_n$ are harmonic in the same neighbourhood of 
$\Gamma$,
the weak convergence implies the uniform convergence of $\nabla v_n$ 
to $\nabla v$ in a neighbourhood of $\Gamma$. As $\nabla v_n=R\,\nabla 
u_n$, the differential forms $-D_2v_n\,dx_1+D_1v_n\,dx_2$ 
are exact, thus 
their integrals over $\Gamma$ vanish. It follows that the integral 
over $\Gamma$ of the differential form $-D_2v\,dx_1+D_1v\,dx_2$ is 
zero. Then we can apply Theorem \ref{due} and we obtain that
the function $u^*$ is a solution to problem (\ref{**}).
 
We now construct a function $w\in \hb(\Omk)$ such that 
$\nabla w=\nabla u^*$ a.e.\ in $\Omega$ and
$w=g$ q.e.\ on $\partial_D\Om\setmeno K$.
As in the proof of Theorem~\ref{uno}, we consider a  connected component 
$C$ of $\Omk$ whose boundary meets $\partial_D\Om\setmeno K$, 
a point $x\in C$, and the corresponding sets $N^{\e}$, $C^{\e}$, 
and ${\Gamma}^{\e}$. Since $K_n\subset N^{\e}$ for $n$ large enough, 
we can apply Corollary~\ref{complete} and deduce that there exists a 
function $w\in\hb(C)$ such that $w=g$ q.e.\ on $\partial 
C\cap(\partial_D\Om\setmeno K)$ and, up to a subsequence,  $\nabla 
u_n\wto \nabla w$ weakly in $L^2(C,\R^2)$, which implies $\nabla 
w=\nabla u^*$ a.e.\ in $C$.

In this way we construct $w$ on all connected components of $\Omk$ 
whose boundary meets $\partial_D\Om\setmeno K$, and we define 
$w:=u^*$ 
on the other connected components. It is clear that $w\in 
\hb(\Omk)$, $\nabla w=\nabla u^*$ a.e.\ in $\Omk$, $w=g$ q.e.\ on 
$\partial_D\Om\setmeno K$, hence 
$w$ is a solution to 
problem (\ref{**})-(\ref{**b}). By the uniqueness of the gradients of 
the solutions, we have $\nabla  u=\nabla w= \nabla u^*$ a.e.\ in $\Omk$.
As the limit does not depend on the subsequence, the whole 
sequence $(\nabla u_n)$ converges to $\nabla u$ weakly in 
$L^2(\Om;\R^2)$.

Taking $u_n-g_n$ and $u-g$ as test functions in the equations satisfied 
by $u_n$ and $u$, we obtain
$$
\int_{\Om}|\nabla u_n|^2dx=\int_{\Om}\nabla u_n \nabla g_n\,dx\,,
\qquad
\int_{\Om}|\nabla u|^2dx=\int_{\Om}\nabla u \,\nabla g\,dx\,.
$$
As $\nabla u_n \wto \nabla u$ weakly in $L^2(\Om,\R^2)$ and $\nabla 
g_n\to\nabla g$ strongly  in $L^2(\Om,\R^2)$, from the 
previous equalities we obtain that $\|\nabla u_n\|_{L^2(\Om,\R^2)}$ 
converges to $\|\nabla u\|_{L^2(\Om,\R^2)}$, which implies the strong 
convergence of the gradients in $L^2(\Om,\R^2)$.
\end{proof}

\end{section}

\begin{section}{HAUSDORFF MEASURE AND HAUSDORFF CONVERGENCE}
\label{Hausdorff}

In this section we study the behaviour of the Hausdorff measure 
$\huno$ along suitable sequences of continua which converge in the 
Hausdorff metric. 
It is well-known that, in general, the Hausdorff measure it is not 
lower semicontinuous with respect to the convergence in the 
Hausdorff metric. When all sets are connected, we have the 
following lower semicontinuity theorem, whose proof can be obtained 
as in Theorem 10.19 of \cite{MS}.

\begin{theorem}[Golab's Theorem]\label{Golab}
Let $(K_n)$ be a sequence of  continua in $\R^2$ which 
converges to a continuum $K$ in the Hausdorff metric. Then
$$
\huno(K\cap U)\le \liminf_{n\to\infty} \,\huno(K_n\cap U)
$$
for every open set $U\subset \R^2$.
\end{theorem}

\begin{remark}\label{connected}{\rm If $C$ is a 
connected subset of $\R^2$ with $\huno(C)<+\infty$, then every 
point in $C$ has a fundamental system of compact  connected 
neighbourhoods (see, e.g., \cite[Lemma 1]{Ch-D}). 
This implies that any two points $x,\,y\in C$ can be 
connected by a continuum contained in $C$. Hence $C$ is the union of an 
increasing sequence $(K_n)$ of continua. As $K_n\to \overline C$ in 
the Hausdorff metric, by Golab's Theorem~\ref{Golab} we have
$\huno(\overline C)\leq\liminf_n\huno(K_n)\leq\huno(C)$. Therefore 
$\huno(\overline C)=\huno(C)$ for every connected subset $C$ of $\R^2$.}
\end{remark}

We shall use also the following consequence of 
Golab's Theorem~\ref{Golab}.
\begin{corollary}\label{sci2}
Let $(H_n)$ be a sequence of compact sets in $\R^2$ which 
converges to a  compact set $H$ in the Hausdorff metric. 
Let  $(K_n)$ be a sequence of continua which 
converges to a continuum $K$ in the Hausdorff metric. Assume that
$H_n\subset K_n$ for every $n$. Then 
\begin{equation}\label{scidif}
\huno(K\setmeno H)\leq
\liminf_{n\to\infty} \, \huno(K_n\setmeno H_n)\,.
\end{equation}
\end{corollary}
\begin{proof}
Given $\e>0$, let 
$H^{\e}:=\{x\in\R^2:{\rm dist}(x,H)\le \e\}$. As 
$H_n\subset  H^{\e}$ for $n$ large enough, we have
$K_n\setmeno H^{\e}\subset  K_n\setmeno H_n$.
Applying Theorem \ref{Golab} with $U=\R^2\setmeno H^{\e}$ we get
$$
\huno(K\setmeno H^{\e}) \le 
\liminf_{n\to\infty} \, \huno(K_n\setmeno H^{\e})
\le \liminf_{n\to\infty} \, \huno(K_n\setmeno H_n)\,.
$$
Passing to the limit as $\e\to0$ we obtain (\ref{scidif}).
\end{proof}

In the next section we shall use the following result.

\begin{lemma}\label{differ}
Let $(H_n)$ be a sequence of continua in $\R^2$
which converges in the 
Hausdorff metric to a  continuum $H$ with $\huno(H)<\pinfty$,
and let $K$ be a  continuum in $\R^2$ with $K\supset H$ and
$\huno(K)<+\infty$.
Then there exists a sequence $(K_n)$ of continua such that
$K_n\to K$ in the Hausdorff metric,
$H_n\subset K_n$, and
$\huno(K_n\setmeno H_n)\to\huno(K\setmeno H)$.
\end{lemma}
\begin{proof} 
It is enough to prove the lemma when $K\setmeno H\neq\emptyset$ 
and $H\neq \emptyset$.
Since $K$ is separable and locally connected (see, e.g.,  
\cite[Lemma 1]{Ch-D}), and 
$K\setmeno H$ is open in $K$, the connected components 
of $K\setmeno H$ are open in $K$ and form a finite or countable 
sequence $(C_i)$. Since each $C_i$ is closed in $K\setmeno H$, we have
$C_i=\overline C_i\cap (K\setmeno H)$.
If $C_i=\overline C_i$, then $K$ would contain an open, closed, and 
non-empty
proper subset, which contradicts the fact that $K$ is connected.
Therefore $C_i\neq \overline C_i$. As 
$\overline C_i\cap (K\setmeno H)=C_i$, 
we conclude that $\emptyset\neq\overline C_i\setmeno C_i\subset H$.
Therefore for every $i$ there exists a 
point $x_i\in \overline C_i\cap H$.
As $H_n\to H$ in the Hausdorff metric, there exists
$x_i^n\in H_n$ such that
$x_i^n\to x_i$. If there are infinitely many connected components $C_i$, 
there exists a sequence of integers $(k_n)$ tending to 
$\infty$ such that
\begin{equation}\label{aax}
\lim_{n\to\infty} \sum_{i=1}^{k_n} |x_i^n - x_i| = 0\,.
\end{equation}
If there are $k<+\infty$ connected components $C_i$, (\ref{aax}) is true with
$k_n=k$ for every $n$.
Let $S_i^n$ denote the closed segment with end-points $x_i^n$ and $x_i$, and 
let
$$
K_n:=H_n\cup\bigcup_{i=1}^{k_n} \overline C_i
\cup\bigcup_{i=1}^{k_n}  S_i^n\,.
$$
Then the sets $K_n$ are continua, contain $H_n$,
and converge to $K$ in the 
Hausdorff metric. As
$\huno(C_i)=\huno(\overline C_i)$ (Remark~\ref{connected}),
we have
$$
\huno(K_n\setmeno H_n)\le \sum_{i=1}^{k_n} \huno(C_i) +
\sum_{i=1}^{k_n} |x_i^n - x_i|\,,
$$
which, together with (\ref{aax}), yields
\begin{equation}\label{aay}
\limsup_{n\to\infty} \,\huno(K_n\setmeno H_n)\le
\huno(K\setmeno H)\,.
\end{equation}
The opposite inequality for the lower limit follows from 
 Corollary~\ref{sci2}.
 \end{proof}

\end{section}

\begin{section}{SOME DIFFERENTIABILITY PROPERTIES}

For every compact set $K$ in $\R^2$ and 
every $g\in\hb(\Omk)$ we define
\begin{equation}\label{e}
\E(g,K):=\min_{v\in {\mathcal V}(g,K)}\Big\{\int_{\Om\setminus K}|\nabla 
v|^2dx+\huno(K)\Big\}\,,
\end{equation}
where ${\mathcal V}(g,K)$ is the set introduced  in (\ref{vgk}).

Given a Hilbert space $X$, we recall that $AC([0,1];X)$  is the space of 
all absolutely continuous functions defined in $[0,1]$ with 
values in $X$. For the main properties of these functions we refer, e.g., 
to \cite[Appendix]{Bre}. Given $g\in AC([0,1];X)$, the time derivative of 
$g$, which exists a.e.\ in $[0,1]$, is denoted by $\gdot$. It is 
well-known that $\gdot$ is a Bochner integrable function with values 
in $X$. 

Finally, let $\K$ be the set of all continua $K$ contained in 
$\overline\Om$ with $\huno(K)<\pinfty$.

The main result of this section is the following theorem.
\begin{theorem}\label{diff5}
Let $g\in AC([0,1];H^1(\Om))$ and let $K\colon[0,1]\to\K$ be an 
increasing function.  Suppose that the function  $t\mapsto \E(g(t),K(t))$ is absolutely 
continuous on $[0,1]$.
Then the  following conditions are equivalent:
\begin{itemize}
\item[(a)] 
\hfil$\displaystyle\frac{d}{ds}\E(g(t),K(s))\Big|\lower1.5ex\hbox{$\scriptstyle s=t$}=0\quad\hbox{ for a.e.\ }
t\in[0,1]\,,$\hfil
\item[(b)] \hfil$\displaystyle\frac{d}{dt}\E(g(t),K(t))=
2(\nabla u(t)|\nabla \gdot(t))\quad \hbox{ for a.e.\ } t\in[0,1]$,\hfil
\end{itemize} 
where $u(t)$ is a solution of the minimum problem (\ref{e}) which defines 
$\E(g(t),K(t))$, and $(\cdot|\cdot)$ denotes the scalar product 
in $L^2(\Om;\R^2)$.
\end{theorem}

To prove Theorem~\ref{diff5} we need the following lemmas.
\begin{lemma}\label{diff1}
Let $K$ be a compact subset of $\overline\Om$ and let 
$ F\colon H^1(\Om)\to\R$ be defined by $ F(g)=\E(g,K)$ for every $g\in 
H^1(\Om)$. Then $ F$ is of class $C^1$ and for every $g,h\in 
H^1(\Om)$ we have 
\begin{equation}\label{dphi}
d F(g)\,h=2\int_{\Om\setminus K}\nabla u_g\nabla h\,dx\,,
\end{equation}
where $u_g$ is a solution of the minimum  problem (\ref{e}) which 
defines $\E(g,K)$.
\end{lemma}
\begin{proof} Since $u_g$ is a solution of problem (\ref{**}) which 
satisfies the boundary condition (\ref{**b}), by linearity 
for every $t\in\R$ we have $\nabla u_{g+th}=\nabla u_g+t\nabla u_h$ 
a.e.\ in $\Om$, hence 
\begin{eqnarray*}
 & F(g+th)- F(g)=\displaystyle \int_{\Om\setminus K}|\nabla u_g+t\nabla 
u_h|^2\,dx-\int_{\Om\setminus K}|\nabla u_g|^2\,dx=\\
&=2t\displaystyle\int_{\Om\setminus K}\nabla u_g\nabla 
u_h\,dx+t^2\int_{\Om\setminus K}|\nabla 
u_h|^2\,dx=2t\int_{\Om\setminus K}\nabla u_g\nabla h\,dx+
t^2\int_{\Om\setminus K}|\nabla u_h|^2\,dx\,,
\end{eqnarray*}
where the last equality is deduced from (\ref{**}). Dividing by $t$  
and letting $t$ tend to $0$ we obtain (\ref{dphi}). The continuity of 
$g\mapsto\nabla u_g$ implies that $ F$ is of class $C^1$.
\end{proof}

Let us consider now the case of time dependent continua $K(t)$.

\begin{lemma}\label{diff2}
Let $K\colon [0,1]\to\K$ be a function and 
let $ F\colon H^1(\Om){\times}[0,1]\to\R$ be defined by 
$ F(g,t)=\E(g,K(t))$. Then the differential $d_1 F$ of $ F$ with 
respect to $g$ is continuous at every 
point  $(g,t)\in H^1(\Om){\times}[0,1]$ such that 
$K(s)\to K(t)$ in the Hausdorff metric as $s\to t$.
\end{lemma}
\begin{proof}
It is enough to apply Lemma~\ref{diff1} and Theorem~\ref{convsol}.
\end{proof}

In the following lemma we find a set of points $t$ satisfying the 
condition considered above.

\begin{lemma}\label{diff3}
Let $K\colon [0,1]\to\K$ be an increasing function and let $t_0$ be a 
continuity point of the 
function $t\mapsto\huno(K(t))$. Then $K(t)\to 
K(t_0)$ in the Hausdorff metric as $t\to t_0$, unless 
$\bigcup_{t<t_0}K(t)=\emptyset\neq \bigcap_{t>t_0}K(t)$. 
\end{lemma}
\begin{proof}
It is enough to show that $K(t_n)\to K(t_0)$ in the Hausdorff metric 
for every monotone sequence $t_n\to t_0$. Let us consider an 
increasing sequence $t_n\nearrow t_0$, and let $K^-$ be the closure 
of the union of the sets $K(t_n)$. Then $K^-$ is a continuum contained 
in $K(t_0)$ and $K(t_n)\to K^-$ in the Hausdorff metric. By the 
continuity of the function $t\mapsto \huno(K(t))$ and by 
Remark~\ref{connected} we have 
$\huno(K^-)=\huno(K(t_0))$, which implies $K^-=K(t_0)$ (see 
Remark~\ref{continua}).
 
Let us consider now a decreasing sequence $t_n\searrow t_0$, and let 
$K^+$ be the intersection of the sets $K(t_n)$. Then $K^+$ is a 
continuum containing $K(t_0)$ and $K(t_n)\to K^+$ in the Hausdorff 
metric. The conclusion can be obtained as in the previous case.
\end{proof}

To deal with the dependence on $t$ of both arguments we need the 
following result.

\begin{lemma}\label{diff4}
Let $X$ be a Hilbert space, let $g\in AC([0,1];X)$, and let 
$ F\colon X{\times}[0,1]\to\R$ be a function such that $ F(\cdot,t)\in 
C^1(X)$ for every $t\in[0,1]$, with differential denoted by 
$d_1 F(\cdot,t)$. Let $t_0\in[0,1]$, let $\psi(t):= F(g(t),t)$, 
and let $\psi_0(t):= F(g(t_0),t)$. Assume that $t_{0}$ is a 
differentiability point of $\psi$ and $g$ and a Lebesgue point of 
$\gdot$, 
 and that $d_1 F$ is continuous at 
$(g(t_0),t_0)$. Then $\psi_0$ is differentiable at $t_0$ and 
$$
\dot\psi_0(t_0)=\dot\psi(t_0)-d_1\! F(g(t_0),t_0)\,\gdot(t_0)\,.
$$
\end{lemma}
\begin{proof}
For every $t\in[0,1]$ we have
\begin{eqnarray*}
\psi_0(t)-\psi_0(t_0)=&\hskip-0.5em F(g(t_0),t)- F(g(t),t)+\psi(t)-\psi(t_0)=\\
=&\hskip-1em \displaystyle\int_t^{t_0}d_1 F(g(s),t)\,\gdot(s)\,ds+\psi(t)-\psi(t_0)\,.
\end{eqnarray*}
The conclusion follows dividing by $t-t_0$ and taking the limit as 
$t\to t_0$.
\end{proof}

\begin{proof}[Proof of Theorem~\ref{diff5}.] 
Let $ F\colon H^1(\Om){\times}[0,1]\to\R$ be defined by 
$ F(g,t)=\E(g,K(t))$. By Lemmas~\ref{diff2} and~\ref{diff3} $d_1 F$ 
is continuous in $(g,t)$ for a.e.\ $t\in[0,1]$ and every $g\in 
H^1(\Om)$. By Lemmas~\ref{diff1} and~\ref{diff4} 
\begin{equation}\label{formula1}
\frac{d}{ds}\E(g(t),K(s))\Big|\lower1.5ex\hbox{$\scriptstyle s=t$}=
\frac{d}{dt}\E(g(t),K(t))-2(\nabla 
u(t)|\nabla \gdot(t))\quad\hbox{ for a.e.\ }t\in[0,1]\,.
\end{equation}
The equivalence between (a) and (b) is now obvious.
\end{proof}

\end{section}

\begin{section}{IRREVERSIBLE QUASI-STATIC EVOLUTION}\label{irrev}

In this section we prove the main result of the paper.
\begin{theorem}\label{kt}
Let $g\in AC([0,1];H^1(\Om))$ and let $K_{0}\in\K$. Then 
there exists a function $K\colon[0,1]\to\K$ such that
\smallskip
\begin{itemize}
\item[(a)] 
\hfil $\displaystyle \vphantom{\frac{d}{ds}} 
K_0\subset K(s)\subset K(t)$  for $0\le s\le t\le 1$, \hfil
\item[(b)] \hfil $\displaystyle \vphantom{\frac{d}{ds}}
\E(g(0),K(0))\leq \E(g(0),K)
\quad\forall\, K\in\K,\,\  K\supset K_0$,\hfil
\item[(c)] \hfil $\displaystyle \vphantom{\frac{d}{ds}}
\hbox{for  }\, 0\le t \le1\quad\E(g(t),K(t))\leq \E(g(t),K)
\quad\forall \, K\in\K,\,\  K\supset K(t)$,\hfil
\item[(d)]\hfil $\displaystyle \vphantom{\frac{d}{ds}}
t\mapsto \E(g(t),K(t)) \hbox{ is 
absolutely continuous on }[0,1]$, \hfil
\item[(e)]\hfil$\displaystyle\frac{d}{ds}\E(g(t),K(s))\Big|\lower1.5ex\hbox{$\scriptstyle s=t$}=0\quad
\hbox{for a.e.\ }t\in[0,1]$.\hfil 
\end{itemize}
\smallskip
Moreover every function $K\colon[0,1]\to\K$ which satisfies (a)--(e) satisfies also
\begin{itemize}
\item[(f)] \hfil$\displaystyle\frac{d}{dt}\E(g(t),K(t))=
2(\nabla u(t)|\nabla \gdot(t))\quad \hbox{ for a.e.\ } t\in[0,1]$,\hfil
\end{itemize}
where $u(t)$ is a solution of the minimum  problem (\ref{e}) which defines 
$\E(g(t),K(t))$.
\end{theorem}
Here and in the rest of the section $(\cdot|\cdot)$  and $\|\cdot\|$  
denote the scalar product  and the norm in $L^2(\Om;\R^2)$.

Theorem~\ref{kt} will be proved by a time discretization process.
Given $\delta>0$, let $N_\delta$ be the largest integer such that
$\delta N_\delta\le 1$; 
for $i\geq 0$ let $t_i^\delta:=i\delta$ and, for $0\le i\le  
N_\delta$, let $g_i^\delta:=g(t_i^\delta)$. 
We define $K_i^\delta$, inductively, as 
a solution of the minimum problem
\begin{equation}\label{pidelta}
\min_K\big\{ \E(g_i^\delta,K) : K\in\K,\ K\supset K_{i-1}^\delta\big\}\,,
\end{equation}
where we set $K_{-1}^\delta:=K_0$.



\begin{lemma}
There exists a solution of the minimum problem (\ref{pidelta}).
\end{lemma}
\begin{proof}
Consider a minimizing sequence $(K_n)$. By the Compactness 
Theorem~\ref{compactness}, passing to a subsequence, we may 
assume that  $(K_n)$ converges in the Hausdorff metric  to
some  continuum  $K$ containing $K_{i-1}^\delta$.
For every $n$ let $u_n$ be a 
solution of the minimum  problem (\ref{e}) which defines $\E(g_i^\delta,K_n)$. 
By Theorem~\ref{convsol} we conclude that $(\nabla u_n)$ converges 
strongly in $L^2(\Om;\R^2)$ to $\nabla u$ where $u$ is a  solution of 
the minimum problem (\ref{e}) which defines $\E(g_i^\delta,K)$.
By Golab's Theorem~\ref{Golab} we have 
$\huno(K)\leq\liminf_{n}\huno(K_n)$. As  
$\|\nabla u\|=\lim_{n}\|\nabla u_n\|$, 
we conclude that $\E(g_i^\delta,K)\leq\liminf_{n}\E(g_i^\delta,K_n)$.  Since 
$(K_n)$ is a minimizing sequence, this proves that $K$  is a 
solution of the minimum problem (\ref{pidelta}). 
\end{proof}

We define now the step functions $g_\delta$, $K_\delta$, and $u_\delta$ 
on $[0,1]$ 
by setting $g_\delta(t):=g_{i}^\delta$, $K_\delta(t):=K_{i}^\delta$, 
and $u_\delta(t):=u_{i}^\delta$ for
$t_i^\delta\leq t<t_{i+1}^\delta$, where $u_i^\delta$ is a  
solution of the minimum problem (\ref{e}) which defines 
$\E(g_{i}^\delta,K_{i}^\delta)$. 

\begin{lemma}\label{discr}
There exists a positive function $\rho(\delta)$, converging to zero 
as $\delta\to0$, such that 
\begin{equation}\label{2discr}
\|\nabla u_j^\delta\|^2+\huno(K_j^\delta)\leq
\|\nabla u_i^\delta\|^2+\huno(K_i^\delta)+
2\int_{t_i^{\delta}}^{t_j^{\delta}}(\nabla u^\delta(t)|\nabla 
\gdot(t))\,dt+\rho(\delta)
\end{equation}
for $0\leq i<j\leq N_\delta$.
\end{lemma}
\begin{proof}
Let us fix an integer $r$ with $i\leq r<j$. {}From the absolute 
continuity of $g$ we have 
$$
g_{r+1}^\delta-g_r^\delta=\int_{t_r^{\delta}}^{t_{r+1}^{\delta}}\gdot(t)\,dt\,,
$$
where the integral is a Bochner integral for functions with values in 
$H^1(\Om)$. This implies that 
\begin{equation}\label{nabla}
\nabla g_{r+1}^\delta-\nabla g_r^\delta=\int_{t_r^{\delta}}^{t_{r+1}^{\delta}}
\nabla \gdot(t)\,dt\,,
\end{equation}
where the integral is a Bochner integral for functions with values in 
$L^2(\Om;\R^2)$.

As $u_r^{\delta}+g_{r+1}^\delta-g_r^\delta\in \hb(\Omk_{r+1}^\delta)$ 
and $u_r^{\delta}+g_{r+1}^\delta-g_r^\delta=g_{r+1}^\delta$ q.e.\ on 
$\partial_D\Omk_r^\delta\supseteq \partial_D\Omk_{r+1}^\delta$, 
from the minimality of $u_{r+1}^\delta$ we obtain,
using (\ref{nabla}), 
\begin{eqnarray*}
& \|\nabla u_{r+1}^\delta\|^2+\huno(K_{r+1}^\delta)\leq
\|\nabla u_r^\delta+\nabla g_{r+1}^\delta-\nabla g_r^\delta\|^2+
\huno(K_r^\delta)\leq\\
&\leq\|\nabla u_{r}^\delta\|^2+\huno(K_r^\delta)+2
\displaystyle\int_{t_r^{\delta}}^{t_{r+1}^{\delta}}(\nabla u_{r}^\delta|\nabla 
\gdot(t))\,dt+\Big(\displaystyle\int_{t_r^{\delta}}^{t_{r+1}^{\delta}}\|\nabla 
\gdot(t)\|\,dt\Big)^2\leq\\
&\leq\|\nabla u_{r}^\delta\|^2+\huno(K_r^\delta)+2
\displaystyle\int_{t_r^{\delta}}^{t_{r+1}^{\delta}}(\nabla u_\delta(t)|\nabla 
\gdot(t))\,dt+\sigma(\delta)\displaystyle\int_{t_r^{\delta}}^{t_{r+1}^{\delta}}\|\nabla 
\gdot(t)\|\,dt\,,
\end{eqnarray*}
where 
$$
\sigma(\delta):=\max_{i\leq r<j} \int_{t_r^{\delta}}^{t_{r+1}^{\delta}}\|\nabla 
\gdot(t)\|\,dt\,.
$$
Iterating now this inequality for $i\leq r<j$ we get 
(\ref{2discr}) with $\rho(\delta):=\sigma(\delta)\int_0^1\|\nabla 
\gdot(t)\|\,dt$.
\end{proof}

\begin{lemma}\label{estim}
There exists a	constant $C$, depending only on $g$ and $K_0$, such that
\begin{equation}\label{stima}																								%
\|\nabla u_i^\delta\|\leq C\quad\hbox{and}\quad\huno(K_i^{\delta})\leq C
\end{equation}
for every $\delta>0$ and for every  $0\leq i\leq N_\delta$.	
\end{lemma}
\begin{proof} 
As $g_i^\delta$ is admissible for the problem 
(\ref{e}) which defines $\E(g_i^\delta, K_i^\delta)$, 
by the minimality of $u_i^\delta$ we have 
$\|\nabla u_i^\delta\|\leq\|\nabla g_i^\delta\|$, hence 
$\|\nabla u_\delta(t)\|\leq\|\nabla g_\delta(t)\|$ for every 
$t\in[0,1]$.
As $t\mapsto g(t)$ is absolutely continuous with values in $H^1(\Om)$ 
the function $t\mapsto \|\nabla\gdot(t)\|$ is integrable on $[0,1]$ 
and 
there exists a constant $C>0$ such that 
$\|\nabla g(t)\|\leq C$ for every $t\in[0,1]$.  This implies the 
former 
inequality in (\ref{stima}). The latter inequality follows 
now from Lemma~\ref{discr}. 
\end{proof}

Let ${S}_\delta$ be the continuum in $\R^3$ defined by
$$
{S}_\delta:=
\bigcup_{0\leq t\leq 1} (\{t\}{\times} K_\delta(t))
=\bigcup_{i=0}^{N_\delta}
([t_i^\delta,\tau_{i+1}^\delta]{\times} K_i^\delta)\,,
$$
where $\tau_{i}^\delta:=t_i^\delta$ for $0\leq i\leq N_\delta$ and 
$\tau_{i}^\delta:=1$ for $i=N_\delta+1$.
By the Compactness Theorem~\ref{compactness} there exists a continuum
${S}\subset[0,1]{\times}\overline\Om$  such that
${S}_{\delta}\to{S}$ in the Hausdorff metric as 
$\delta\to 0$ along a suitable sequence. 
For every $t\in[0,1]$ we define 
\begin{equation}\label{k^*}
\textstyle K^+(t):=\{x\in\R^2:(t,x)\in{S}\}\,.
\end{equation}
As $K_\delta(t)$ is increasing with respect to $t$, the same is true 
for $K^+(t)$.
For $0<t\le 1$ we define 
\begin{equation}\label{k_*}
\textstyle K^-(t):={\rm cl}\big(\bigcup_{s<t}K^+(s)\big)\,,
\end{equation}
where ${\rm cl}$ denotes the closure.
By the monotonicity of $K^+(t)$ and $K^-(t)$, there is at most 
one point $t_0\in(0,1]$ such that
$K^-(t_0)=\emptyset$ and $K^+(t_0)\neq \emptyset$.
Let 
$\Theta$ be the set of points $t\in(0,1)$ such that  $K^-(t)=K^+(t)$,
and let 
$\Theta^*$ be the set of continuity points in $(0,1)$ of the function 
$t\mapsto \huno(K^+(t))$.

For the reader's convenience we give here the complete proof of the 
following result, which was originally proved in  \cite[Proposition 4]{ChD}. 
\begin{lemma}\label{ChD} For every
$t\in[0,1]$ the set $K^+(t)$ is a continuum and 
$K^+(t)= \bigcap_{s>t}K^+(s)$ for $0\le t<1$.
For every $t\in \Theta$ the sequence $K_{\delta}(t)$ 
converges to $K^+(t)$ in the 
Hausdorff metric as $\delta\to0$ along the same sequence considered 
for ${S}_\delta$.
Moreover $[0,1]\setmeno\Theta$ is at most countable and
$\Theta=\Theta^*\setmeno\{t_0\}$.
\end{lemma}
\begin{proof} 
By construction $K^+(t)$ is compact. To prove that it is connected, let 
$x,\, y\in K^+(t)$. Then $(t,x),\, (t,y)\in {S}$ and, since ${S}$ is the  
limit of ${S}_{\delta}$ in the Hausdorff metric 
as $\delta\to0$
along a suitable sequence, there 
exist sequences $(s_\delta, x_\delta)$, 
$(t_\delta,y_\delta)\in {S}_{\delta}$ converging 
to $(t,x)$ and $(t,y)$, respectively. By the monotonicity of 
$s\mapsto K_{\delta}(s)$,
we may assume that
$s_\delta=t_\delta$, so that $x_\delta$ and $y_\delta$ belong to the same 
continuum 
$K_{\delta}(t_\delta)$. Passing to a subsequence, we
may assume that $K_{\delta}(t_\delta)$ converges in the Hausdorff metric
to some continuum $K$. Then $\{t_\delta\}{\times} K_{\delta}(t_\delta)$ 
converges to $\{t\}\times K$, which implies that 
$\{t\}{\times}K\subset{S}$, so
that $K\subset K^+(t)$. As $x$ and $y$ are connected in $K^+(t)$ by the set
$K$, it follows that $K^+(t)$ is connected.

The monotonicity of $s\mapsto K^+(s)$ implies that 
$K^+(t)\subset\bigcap_{s>t}K^+(s)$ for $0\le t<1$.
To prove the reverse inclusion, let now $x\in\bigcap_{s>t}K^+(s)$ and
$t_n\searrow t$. Then $(t_n,x)\in {S}$ and, as
${S}$ is closed, $(t,x)\in{S}$, which implies that 
$x\in K^+(t)$.

Let us fix $t\in \Theta$. To prove that  $K_{\delta}(t)$ 
converges to $K^+(t)$, by the Compactness 
Theorem~\ref{compactness} we may assume that $K_{\delta}(t)$ 
converges  in the Hausdorff metric to some compact set $K$ as 
$\delta\to0$ along a 
suitable subsequence of the sequence considered for ${S}_\delta$, 
and we have only to prove that 
$K^-(t)\subset K\subset K^+(t)$. 

Let us prove the first inclusion. As $K$ is closed, it is enough to show that 
$\bigcup_{s<t}K^+(s)\subset K$. Let $x\in \bigcup_{s<t}K^+(s)$. Then there 
exists $ s<t$ such that  $(s,x)\in {S}$. As 
${S}_{\delta}\to {S}$ in the Hausdorff metric, 
there exists a sequence $(s_\delta,x_\delta)\to(s,x)$ such that 
$(s_\delta,x_\delta)\in{S}_{\delta}$, and hence $x_\delta\in 
K_{\delta}(s_\delta)$. As $s_\delta<t$ for $\delta$ small enough, by the 
monotonicity of $s\mapsto K_\delta(s)$, we have $x_\delta\in 
K_{\delta}(t)$. Since $x_\delta\to x$ we conclude that $x\in K$. 

Let us prove that $K\subset K^+(t)$. Let $x\in K$. By the definition 
of $K$ there exist $x_\delta\in K_{\delta}(t)$ such that $x_\delta\to x$. 
Then $(t,x_\delta)\in {S}_{\delta}$ and $(t,x_\delta)\to(t,x)$. This 
implies $(t,x)\in {S}$, hence $x\in K^+(t)$. 
 
It is easy to see that for $0<t<1$
\begin{eqnarray*}
&\displaystyle \lim_{s\to t+}\huno(K^+(s))=\huno(K^+(t))\,,\\
& {\displaystyle \lim_{s\to 
t-}\huno(K^+(s))}=\huno(\textstyle\bigcup_{s<t}K^+(s))=\huno(K^-(t))\,,
\end{eqnarray*}
where the last equality follows from Remark~\ref{connected}. This implies 
that $t\in\Theta^*$ if and only if $\huno(K^+(t)\setmeno K^-(t))=0$.
The equality $\Theta=\Theta^*\setmeno\{t_0\}$ follows 
now from Remark~\ref{continua}.
The statement about $[0,1]\setmeno\Theta$ is a consequence of this 
equality and of the  
continuity properties of monotone functions.
\end{proof}

Since $[0,1]\setmeno\Theta$  is at most countable, by a diagonal 
argument we may extract a subsequence of the sequence considered for
${S}_\delta$  such that, for every $t\in[0,1]\setmeno\Theta$, 
$K_{\delta}(t)$ converge to a continuum 
$K(t)$ in 
the Hausdorff metric as $\delta\to0$ along this subsequence. In the 
rest of this section we shall always refer to this subsequence when 
we write $\delta\to0$. 
Let us define $K(t):=K^+(t)$ for all $t\in \Theta$. By the previous 
discussion and by Lemma~\ref{ChD} for every $t\in[0,1]$ we have  
$K_{\delta}(t)\to K(t)$ in the Hausdorff metric as $\delta\to0$, 
hence $K(t)$ is a continuum.  
By Lemma~\ref{estim} there exist a constant $C$ such that 
$\huno(K_{\delta}(t))\leq C$ for every $t$. By 
Golab's Theorem \ref{Golab} we have $\huno(K(t))\leq C$. This shows 
that $K(t)\in \K$ for every $t\in[0,1]$.

Since $K_0\subset K^0_{\delta}=K_{\delta}(0)$, we have that 
$K_0\subset K(0)$.
Since $K_{\delta}(t)$ is increasing with respect to $t$, the same 
property holds for $K(t)$.
This shows 
that $K(t)$ satisfies condition (a) of Theorem~\ref{kt}. 

Moreover it is easy to see that
\begin{eqnarray}
&\textstyle K^+(t)=\bigcap_{s>t}K(s)\qquad\hbox{for }0
\leq t<1\,,\label{k^*t}\\
&\textstyle K^-(t)={\rm cl}\big(\bigcup_{s<t}K(s)\big)
\qquad\hbox{for }0<t\leq1\,,
\label{k_*t}\\
&\textstyle K^-(t)\subset K(t)\subset K^+(t)
\qquad\hbox{for }0< t<1\,.\label{7.9}
\end{eqnarray}
Indeed, (\ref{k^*t}) and (\ref{k_*t}) follow from (\ref{7.9}), which 
has been proved in the previous discussion.

For every $t\in[0,1]$ let $u(t)$ be  a solution of  the minimum problem 
(\ref{e}) which defines $\E(g(t),K(t))$.
\begin{lemma}\label{l2}
For every $t\in[0,1]$ we have 
$\nabla u_{\delta}(t)\to\nabla u(t)$ strongly in $L^2(\Om;\R^2)$.
\end{lemma}
\begin{proof} As $u_{\delta}(t)$ is a solution of the minimum  problem 
(\ref{e}) which defines $\E(g_{\delta}(t),K_{\delta}(t))$, and 
$g_\delta(t)\to g(t)$ strongly in $H^1(\Om)$, the 
conclusion follows from 
Theorem~\ref{convsol}. 
\end{proof}

\begin{lemma}\label{condb}For every $t\in[0,1]$
we have 
\begin{equation}\label{pt}
\E(g(t),K(t))\leq \E(g(t),K)\quad\forall\,K\in\K\,\ K\supset 
K(t)\,.
\end{equation}
Moreover
\begin{equation}\label{pt0}
\E(g(0),K(0))\leq \E(g(0),K)\quad\forall\,K\in\K\,\ K\supset 
K_0\,.
\end{equation}
\end{lemma}
\begin{proof} Let us fix $t\in[0,1]$ and $K\in\K$ with $K\supset K(t)$. 
Since  $K_\delta(t)$ converges
to $K(t)$ in the Hausdorff metric as $\delta\to0$, 
by Lemma~\ref{differ} there exists a sequence of continua 
$(K_\delta)$, 
converging to $K$ in the Hausdorff metric, such that
$K_\delta\supset K_\delta(t)$ and 
$\huno(K_\delta\setmeno K_\delta(t))\to \huno(K\setmeno K(t))$
as $\delta\to0$.

Let $v_\delta$ and $v$ be solutions of the minimum problems 
(\ref{e}) which define $\E(g_\delta(t),K_\delta)$ and 
$\E(g(t),K)$, respectively. By Theorem~\ref{convsol}
$\nabla v_\delta\to\nabla v$ strongly
in $L^2(\Om;\R^2)$. 
The minimality of $K_\delta(t)$ expressed by (\ref{pidelta})
gives 
$\E(g_\delta(t),K_\delta(t))\leq \E(g_\delta(t),K_\delta)$, 
which implies 
$\|\nabla u_\delta(t)\|^2\leq\|\nabla v_\delta\|^2+\huno(K_\delta\setmeno 
K_\delta(t))$.
Passing to the limit as $\delta\to0$ and using Lemma~\ref{l2}
we get $\|\nabla u(t)\|^2\leq\|\nabla v\|^2+\huno(K\setmeno K(t))$.
Adding $\huno(K(t))$ to both sides we obtain (\ref{pt}).

A similar proof holds for (\ref{pt0}).
The only difference is that now we take $K_\delta=K$ and  use the fact that
$\E(g_\delta(0),K_\delta(0)) \le \E(g_\delta(0),K)=\E(g(0),K)$, which implies
$\|\nabla u_\delta(0)\|^2+\huno(K_\delta(0))\leq
\|\nabla v\|^2+\huno(K)$.
Passing to the limit as $\delta\to0$ and using Lemma~\ref{l2}
and Golab's Theorem \ref{Golab} we obtain (\ref{pt0}).
\end{proof}

The previous lemma proves conditions (b) and (c) of Theorem~\ref{kt}. 
To show that  conditions (d) and (e) are also satisfied, 
we begin by proving the following inequality.  
\begin{lemma}\label{ineq}
For every $s,\, t$ with $0\leq s<t\leq1$
\begin{equation}\label{diseg}
\|\nabla u(t)\|^2+\huno(K(t))\leq\|\nabla u(s)\|^2+\huno(K(s))+
2\int_s^t(\nabla u(\tau)|\nabla \gdot(\tau))d\tau\,.
\end{equation}
\end{lemma}
\begin{proof}
Let us fix $s,t$ with $0\leq s<t\leq1$. Given $\delta>0$ let $i$ and $j$ 
be the 
integers such that $t_i^\delta\leq s<t_{i+1}^\delta$ and 
$t_j^\delta\leq t<t_{j+1}^\delta$. 
Let us define $s_\delta:=t_i^\delta$ and $t_\delta:=t_j^\delta$. 
Applying Lemma~\ref{discr} we obtain
\begin{equation}\label{ediscr}
\|\nabla u_\delta(t)\|^2+\huno(K_\delta(t)\setmeno K_\delta(s))\leq
\|\nabla u_\delta(s)\|^2+
2\int_{s_{\delta}}^{t_{\delta}}\!\!\!(\nabla u^\delta(\tau)|\nabla 
\gdot(\tau))\,d\tau+\rho(\delta)\,,
\end{equation}
with $\rho(\delta)$ converging to zero as $\delta\to0$.											  %

Since $K_\delta(\tau)$ 
converges to $K(\tau)$ in the Hausdorff metric,  by Lemma~\ref{l2} 
for every $\tau\in[0,1]$ we have $\nabla u_{\delta}(\tau)\to \nabla u(\tau)$ 
 strongly in 
$L^2(\Om,\R^2)$ as $\delta\to0$.
By Corollary~\ref{sci2} we get
$$
\huno(K(t)\setmeno K(s))\leq\liminf_{\delta\to0}\huno(K_\delta(t)\setmeno 
K_\delta(s))\,.
$$
Passing now to the limit in (\ref{ediscr}) as $\delta\to0$ we 
obtain (\ref{diseg}).
\end{proof}

The following lemma concludes the proof of Theorem~\ref{kt},
showing that also conditions (d), (e) and (f) are satisfied. 

\begin{lemma}
The function $t\mapsto \E(g(t),K(t))$ is absolutely continuous on 
$[0,1]$ and 
\begin{equation}\label{gtkt}
\frac{d}{dt}\E(g(t),K(t))=2(\nabla u(t)|\nabla\gdot(t))\qquad\hbox{ for 
a.e.\ }t\in[0,1]\,.
\end{equation}
Moreover
\begin{equation}\label{gtks}
\frac{d}{ds}\E(g(t),K(s))\Big|\lower1.5ex\hbox{$\scriptstyle s=t$}=0\qquad\hbox{ for a.e.\ }t\in[0,1]\,.
\end{equation}
\end{lemma}
\begin{proof}Let $0\leq s<t\leq1$.
{}From the previous lemma we get 
\begin{equation}\label{alto}
\E(g(t),K(t))-\E(g(s),K(s))\leq2\int_s^t(\nabla 
u(\tau)|\nabla\gdot(\tau))\,d\tau\,.
\end{equation}
On the other hand, by condition (c) of Theorem~\ref{kt} we have 
$\E(g(s),K(s))\leq \E(g(s),K(t))$, and by Lemma~\ref{diff1} 
$$
\E(g(t),K(t))-\E(g(s),K(t))=
2\int_s^t(\nabla u(\tau, t)|\nabla \gdot(\tau))\,d\tau\,,
$$
where $u(\tau,t)$ is a solution of the minimum problem~(\ref{e}) which defines 
$\E(g(\tau),K(t))$.
Therefore
\begin{equation}\label{basso}
\E(g(t),K(t))-\E(g(s),K(s))\geq 
2\int_s^t(\nabla u(\tau, t)|\nabla \gdot(\tau))\,d\tau\,.
\end{equation}
Since there exists a constant $C$ such that
$\|\nabla u(\tau)\|\leq\|\nabla g(\tau)\|\leq C$ and  
 $\|\nabla u(\tau,t)\|\leq\|\nabla g(\tau)\|\leq C$ for $s\leq\tau\leq 
 t$, from (\ref{alto}) and (\ref{basso})
we obtain
$$
\big|\E(g(t),K(t))-\E(g(s),K(s))\big|\leq 2\,C
\int_s^t\|\nabla \gdot(\tau)\|\,d\tau\,,
$$
which proves that the  function $t\mapsto \E(g(t),K(t))$ is absolutely 
continuous.

As $\nabla u(\tau,t)\to\nabla u(t)$ strongly in $L^2(\Om,\R^2)$ when  
$\tau\to t$, if we divide (\ref{alto}) and (\ref{basso}) by $t-s$, and 
take the limit as $s\to t-$ we obtain (\ref{gtkt}). 
Equality (\ref{gtks}) follows from Theorem~\ref{diff5}.
\end{proof}

Theorem~\ref{kt0} is a consequence of Theorem~\ref{kt} and of
the following lemma.
\begin{lemma}\label{previous}
Let $K\colon[0,1]\to\K$ be a function which satisfies 
conditions (a)--(e) of Theorem~\ref{kt}. Then, for $0<t\le1$, 
\begin{equation}\label{s<t}
\E(g(t),K(t))\leq \E(g(t),K)
\quad\forall \,K\in\K\,\ K\supset \textstyle
\bigcup_{s<t}K(s)\,.
\end{equation}
\end{lemma}

\begin{proof}
Let us fix $t$, with $0<t\le1$, and a continuum $K\supset 
\bigcup_{s<t}K(s)$. For  $0\leq s<t$ we have $K\supset K(s)$, and from
condition (c) of Theorem~\ref{kt} we obtain $\E(g(s),K(s))\leq 
\E(g(s),K)$. As the functions $s\mapsto \E(g(s),K(s))$ and
$s\mapsto \E(g(s),K)$ are continuous, passing to the limit as $s\to 
t-$ we get (\ref{s<t}).
\end{proof}

The following lemma shows that $K(t)$, $K^-(t)$, and $K^+(t)$ have the 
same total energy.
\begin{lemma}\label{*****}
Let $K\colon[0,1]\to\K$ be a function which satisfies 
conditions (a)--(e) of Theorem~\ref{kt}, and let $K^-(t)$ and 
$K^+(t)$ be defined by (\ref{k_*t}) and (\ref{k^*t}). Then
\begin{eqnarray}
&&\E(g(t),K(t))=\E(g(t),K^-(t))\quad\hbox{for }0<t\le1\,,\label{eq*}\\
&&\E(g(t),K(t))= \E(g(t),K^+(t))\quad\hbox{for }0\le t<1\,.\label{eq**}
\end{eqnarray}
\end{lemma}

\begin{proof} Let $0<t\le1$. Since $K(s)\to K^-(t)$ in the Hausdorff 
metric as $s\to t-$, and $\huno(K(s))\to \huno(K^-(t))$ by 
Remark~\ref{connected}, it follows that $\E(g(s),K(s))\to \E(g(t),K^-(t))$ as 
$s\to 
t-$  by Theorem~\ref{convsol}. As the function $s\mapsto \E(g(s),K(s))$ 
is continuous, we obtain (\ref{eq*}). The proof of (\ref{eq**}) is 
analogous.
\end{proof}

\begin{remark}\label{rem}
{\rm {}From Lemmas~\ref{previous} and~\ref{*****} it 
follows that, if $K\colon[0,1]\to\K$ is a function which satisfies 
conditions (a)--(e) of Theorem~\ref{kt},  the same is true for the 
functions 
$$
t\mapsto\left\{\begin{array}{ll}
\hskip-0.5em K(0) & \hbox{for }t=0\,,\\
\hskip-0.5em K^-(t)& \hbox{for } 0<t\leq 1
\,,
\end{array}\right.\qquad
t\mapsto\left\{\begin{array}{ll}
\hskip-0.5em K^+(t)& \hbox{for }0\leq t<1\,,\\
\hskip-0.5em K(1)& \hbox{for } t=1
\,,
\end{array}\right.
$$
where $K^-(t)$ and 
$K^+(t)$ are defined by (\ref{k_*t}) and (\ref{k^*t}). Therefore the 
problem has a left-continuous  solution and a 
right-continuous solution.
}\end{remark}

\begin{remark}\label{g0}
{\rm In Theorem~\ref{kt} suppose that
$\E(g(0),K_0)\le \E(g(0),K)$
for every $K\in\K$ with $K\supset K_0$.
Then in our time discretization process we can take
$K_0^\delta=K_0$ for every $\delta>0$.
Therefore there exists a function $K\colon[0,1]\to\K$, 
satisfying conditions (a)--(e) of Theorem~\ref{kt}, such that
$K(0)=K_0$. In particular this happens for every  $K_0$ 
whenever $g(0)=0$.

If $\E(g(0),K_0)< \E(g(0),K)$
for every $K\in\K$ with $K\supset K_0$ and  $K\neq K_0$,
then every function $K\colon[0,1]\to\K$, which satisfies
conditions (a)--(e) of Theorem~\ref{kt}, satisfies also $K(0)=K_0$ by~(b).
This always happens when $g(0)=0$ and $K_0\neq\emptyset$.

In the case $g(0)=0$ and $K_0=\emptyset$,  by condition (b) we must 
have $\huno(K(0))=0$, hence $K(0)$ has at most one element. Besides 
the solution with $K(0)=\emptyset$, considered in the first part of 
this remark, for every 
$x_0\in\overline\Om$ there exists a solution with 
$K(0)=\{x_0\}$. Indeed in our time discretization process we can take
$K_0^\delta=\{x_0\}$ for every $\delta>0$.
}\end{remark}

We consider now the case where $g(t)$ is proportional to a fixed 
function $h\in H^1(\Om)$.

\begin{proposition}\label{th}
In Theorem~\ref{kt} suppose that $g(t)=\varphi(t)\,h$, where 
$\varphi\in AC([0,1])$ is non-decreasing and non-negative, and $h$ is 
a fixed function in $H^1(\Om)$. Let $K\colon[0,1]\to\K$ be a function which satisfies 
conditions (a)--(e) of Theorem~\ref{kt}. Then
\begin{equation}
\E(g(t),K(t))\le \E(g(t),K(s))
\end{equation}
for $0\le s<t\le1$.
\end{proposition}

\begin{proof} Let us fix $0\le s<t\le1$. For every $\tau\in[0,1]$ let 
$v(\tau)$ be a solution of the minimum problem (\ref{e}) which defines 
$\E(h,K(\tau))$. As $u(\tau)=\varphi(\tau)\,v(\tau)$ and 
$\gdot(\tau)=\dot\varphi(\tau)\,h$, 
from condition (f) we obtain, adding and subtracting
$\E(g(s),K(s))$,
\begin{eqnarray*}
&\E(g(t),K(t))-\E(g(t),K(s))=\\
&\displaystyle=2\int_s^t(\nabla v(\tau)|\nabla h)
\,\varphi(\tau)\,\dot{\varphi}(\tau)\,d\tau +
(\varphi(s)^2-\varphi(t)^2)\|\nabla v(s)\|^2\,.
\end{eqnarray*}
As $v(\tau)$ is a solution of problem (\ref{**})
with $K=K(\tau)$, and
$v(\tau)=h$ q.e.\ on $\partial_D\Omk(\tau)$, we have
$(\nabla v(\tau)|\nabla h)=\|\nabla v(\tau)\|^2$.
By the monotonicity of $\tau\mapsto K(\tau)$, for $s\le \tau\le t$ 
we have $v(s)\in \hb(\Omk(\tau))$ and $v(s)=h$ 
q.e.\ on $\partial_D\Omk(\tau)$. By the minimum property of $v(\tau)$
we obtain $\|\nabla v(\tau)\|^2\le \|\nabla v(s)\|^2$ for $s\le 
\tau\le t$. 
Therefore
\begin{eqnarray*}
&\E(g(t),K(t))-\E(g(t),K(s))\le\\
&\displaystyle
\le 2 \int_s^t \varphi(\tau)\,\dot{\varphi}(\tau)\, d\tau \;
\|\nabla v(s)\|^2 +
(\varphi(s)^2-\varphi(t)^2)\|\nabla v(s)\|^2= 0\,,
\end{eqnarray*}
which concludes the proof.
\end{proof}
\end{section}

\begin{section}{BEHAVIOUR NEAR THE TIPS}\label{tips}

In this section we consider a function $K\colon[0,1]\to\K$ 
which satisfies conditions (a)--(e) of Theorem~\ref{kt} for a suitable $g\in 
AC([0,1];H^1(\Om))$, and we study the behaviour of 
the solutions $u(t)$ near the ``tips'' of the sets $K(t)$. Under 
some natural assumptions, we shall see that $K(t)$ satisfies Griffith's 
criterion for crack growth.

For every bounded open set $A\subset\R^2$ with Lipschitz boundary, 
for every compact set $K\subset\R^2$, and for every function 
$g\colon \partial A\setmeno K\to\R$ we define
\begin{equation}\label{eloc}
\E(g,K, A):=\min_{v\in{\mathcal V}(g,K,A)}
\Big\{\int_{A\setminus K}|\nabla v|^2\,dx
+\huno(K\cap \overline A)\Big\} \,,
\end{equation}
where 
$$
{\mathcal V}(g,K,A):=\{v\in\hb(A\setmeno K): v=g \quad\hbox{q.e.\ on }
\partial A\setmeno K\}\,.
$$

We now consider in particular the case where $K$ is a regular arc, and
summarize some known results on the behaviour of a solution of 
problem (\ref{**}) near the end-points of $K$. Let $B$ be an open ball 
in $\R^2$ and let $\gamma\colon [\sigma_0, \sigma_1]\to \R^2$ be a 
simple path of class $C^2$ parametrized by arc length. Assume that 
$\gamma(\sigma_0)\in\partial B$ and $\gamma(\sigma_1)\in\partial B$, 
while $\gamma(\sigma)\in B$ for $\sigma_0< \sigma< \sigma_1$. Assume 
in addition that $\gamma$ is not tangent to $\partial B$ at $\sigma_0$
and $\sigma_1$.
For every $\sigma\in [\sigma_0, \sigma_1]$ let ${\Gamma}(\sigma):=
\{\gamma(s): \sigma_0\le s \le \sigma\}$.

\begin{theorem}\label{Grisvard1}
Let $\sigma_0< \sigma< \sigma_1$ and let $u$ be a solution to 
problem (\ref{**}) with $\Omega=B$,  $\partial_D\Omega=\partial B$, 
and $K={\Gamma}(\sigma)$.
Then there exists a unique constant $\kappa=\kappa(u,\sigma)\in\R$ such that
\begin{equation}\label{sif}
u-\kappa \,\sqrt{2\rho/\pi}\,
\sin(\theta/2)\in H^2(B\setmeno {\Gamma}(\sigma)) \cap 
H^{1,\infty}(B\setmeno {\Gamma}(\sigma)) \,,
\end{equation}
where $\rho(x)=|x-\gamma(\sigma)|$ and $\theta(x)$ is the continuous function on 
$B\setmeno {\Gamma}(\sigma)$ which coincides with the oriented angle 
between $
\dot\gamma(\sigma)$ and $x-\gamma(\sigma)$, and 
vanishes on the points of the form 
$x=\gamma(\sigma)+\e\, 
\dot \gamma(\sigma)$ for 
sufficiently small $\e>0$.
\end{theorem}

\begin{proof} 
Let $B^-$ and $B^+$ be the connected components of 
$B\setmeno {\Gamma}(\sigma_1)$. Since $B^-$ and $B^+$ have a Lipschitz
boundary, by Proposition~\ref{h1} $u$ belongs to $H^1(B^-)$ and
$H^1(B^+)$. This implies that $u\in L^2(B)$, and hence 
$u\in H^1(B\setmeno {\Gamma}(\sigma))$. The conclusion follows now from
\cite[Theorem 4.4.3.7 and Section 5.2]{Gri1}, as shown in 
\cite[Appendix1]{MSh}.
\end{proof}

\begin{remark}\label{stress}{\rm
If $u$ is interpreted as the third component of the displacement in an 
anti-plane shear, as we did in the introduction, then $\kappa$ 
coincides with the {\it Mode III stress intensity factor\/} $K_{I\!I\!I}$
of the displacement $(0,0,u)$.
}\end{remark}

\begin{theorem}\label{Grisvard2}
Let $g\colon \partial B\setmeno \{\gamma(\sigma_0)\} \to\R$ be a function such that for 
every $\sigma_0< \sigma< \sigma_1$ there exists
$g(\sigma)\in \hb(B\setmeno {\Gamma}(\sigma))$ with $g(\sigma)=g$
q.e.\ on $\partial B\setmeno {\Gamma}(\sigma)=
\partial B\setmeno \{\gamma(\sigma_0)\}$. Let $v(\sigma)$ be a 
solution of the minimum problem (\ref{eloc}) which defines 
$\E(g, {\Gamma}(\sigma), B)$.
Then, for every $\sigma_0< \sigma< \sigma_1$, 
$$
\frac{d}{d\sigma}\E(g, {\Gamma}(\sigma), B)=
1-\kappa(v(\sigma),\sigma)^2\,,
$$
where $\kappa$ is defined by (\ref{sif}).
\end{theorem}

\begin{proof} 
It is enough to adapt the proof of \cite[Theorem 6.4.1]{Gri2}.
\end{proof}

Let us return to the function $K\colon[0,1]\to\K$ 
considered at the beginning of the section, and let 
$0\le t_0<t_1\le 1$. Suppose that the following structure condition is 
satisfied: there exist a finite family of simple arcs ${\Gamma}_i$,
$i=1,\ldots,m$, contained in $\Om$ and parametrized by arc length
by $C^2$  paths $\gamma_i\colon [\sigma_i^0, \sigma_i^1]\to \Om$, 
such that, for $t_0<t<t_1$,
\begin{equation}\label{structure}
K(t)=K(t_0)\cup \bigcup_{i=1}^m {\Gamma}_i(\sigma_i(t))\,,
\end{equation}
where ${\Gamma}_i(\sigma):=
\{\gamma_i(\tau): \sigma_i^0\le \tau \le \sigma\}$ and
$\sigma_i\colon [t_0,t_1]\to [\sigma_i^0, \sigma_i^1]$ are 
non-decreasing functions with $\sigma_i(t_0)=\sigma_i^0$ and
$\sigma_i^0<\sigma_i(t)<\sigma_i^1$ for $t_0<t<t_1$. 
Assume also that the arcs
${\Gamma}_i$ are pairwise disjoint, and that 
${\Gamma}_i\cap K(t_0)=\{\gamma_i(\sigma_i^0)\}$. 
We consider the sets ${\Gamma}_i(\sigma_i(t))$ as the increasing 
branches of the fracture $K(t)$ and the
points $\gamma_i(\sigma_i(t))$ as their moving tips.
For $i=1,\ldots,m$ and 
$\sigma_i^0<\sigma<\sigma_i^1$
let $\kappa_i(u,\sigma)$ be the stress intensity factor
defined by (\ref{sif}) with $\gamma=\gamma_i$ and $B$ equal to a 
sufficiently small 
ball centred at $\gamma_i(\sigma)$.

We are now in a position to state the main result of this section.

\begin{theorem}\label{Griffith}
Let $K\colon[0,1]\to\K$ be a function which satisfies 
conditions (a)--(e) of Theorem~\ref{kt} for a suitable 
$g\in AC([0,1];H^1(\Om))$, let $u(t)$ be a solution of the minimum
problem (\ref{e}) which defines $\E(g(t),K(t))$,
and let $0\le t_0<t_1\le 1$.
Assume that (\ref{structure}) is satisfied for $t_0<t<t_1$, and that 
the arcs ${\Gamma}_i$ and the functions $\sigma_i$
satisfy all properties considered above.
Then
\begin{eqnarray}
&\dot\sigma_i(t)\ge 0\quad \hbox{ for a.e.\ } t\in (t_0,t_1)\,,
\label{sigmadot}\\
&1-\kappa_i(u(t),\sigma_i(t))^2\ge 0
\quad \hbox{ for every\ } t\in (t_0,t_1)\,,
\label{sif>}\\
&\big\{1- \kappa_i(u(t),\sigma_i(t))^2\big\}
\,\dot\sigma_i(t)= 0\quad \hbox{ for a.e.\ } t\in (t_0,t_1)\,,
\label{sif=}
\end{eqnarray}
for $i=1,\ldots,m$.
\end{theorem}

The first condition says simply that the length of every branch
of the fracture can
not decrease, and reflects the irreversibility of the process. 
The second condition says 
that the absolute value of the stress intensity factor must be less than or equal to 
$1$ at every tip and for every time. The last condition says that, 
at a given tip, the stress intensity factor must be equal to $\pm 1$ 
at almost every time in which this tip moves with a positive velocity.
This is Griffith's criterion for crack growth in our model.

To prove Theorem~\ref{Griffith} we use the following lemma.

\begin{lemma}\label{cond3loc}
Let $g\in H^1(\Om)$, let $H\in\K$, let $u$ be a solution of the 
minimum problem (\ref{e}) which 
defines $\E(g,H)$, and let $A$ be an open subset of $\Om$, with 
Lipschitz boundary, such that $H\cap\overline A\neq\emptyset$. Assume that 
\begin{equation}\label{minn}
\E(g,H)\leq \E(g,K)\qquad\forall \,K\in\K,\,\ K\supset 
H\,.
\end{equation}
Then 
\begin{equation}\label{minnA}
\E(u,H,A)\leq \E(u,K,A)\qquad\forall \,K\in\KA,\,\ K\supset 
H\cap\overline A\,,
\end{equation}
where $\KA$ is the set of all continua $K\subset\overline A$ with 
$\huno(K)<\pinfty$.
\end{lemma}

\begin{proof}
Let $K\in\KA$ with $K\supset 
H\cap\overline A$,  let $v$ be  a solution of the 
minimum problem (\ref{eloc})  which defines $\E(u,K,A)$, 
and let $w$ be the function defined 
by $w:=v$ on $\overline A\setmeno K$ and by $w:=u$ on 
$(\overline\Om\setmeno \overline A)\setmeno  H$. As $v=u$ q.e.\ on 
$\partial A\setmeno K$ the function $w$ belongs to $\hb(\Om\setmeno 
(H\cup K))$; using also the fact that $u=g$ q.e.\ on  $\partial_D\Om\setmeno 
H$, we obtain that  $w=g$ q.e.\ on $\partial_D\Om\setmeno (H\cup K)$.
Therefore 
\begin{eqnarray}
&\displaystyle
\E(g,H\cup K)\leq\int_{\Om\setminus (H\cup K)}|\nabla 
w|^2\,dx+\huno(H\cup K)=\label{pprima}\\
&\displaystyle
= \int_{A\setminus K}|\nabla v|^2\,dx+\huno(K\cap \overline A)+
\int_{(\Om\setminus A)\setminus H}|\nabla u|^2\,dx+
\huno(H\setmeno\overline A)\,.\nonumber
\end{eqnarray}
On the other hand, by the minimality of $u$, 
\begin{eqnarray}
&\displaystyle
\int_{A\setminus H}|\nabla u|^2\,dx+\huno(H\cap \overline A)+
\int_{(\Om\setminus A)\setminus H}|\nabla u|^2\,dx+\huno(H\setmeno\overline 
A)=\label{sseconda}\\
&\displaystyle
=\int_{\Om\setminus H}|\nabla u|^2\,dx+\huno(H)=\E(g,H)\leq\E(g,H\cup K)\,,
\nonumber
\end{eqnarray}
where the last inequality follows from (\ref{minn}), since $H\cup K$ 
is connected.
{}From (\ref{pprima}) and (\ref{sseconda}) we obtain 
$$
\int_{A\setminus H}|\nabla u|^2\,dx+ \huno(H\cap \overline A)\leq
\int_{A\setminus K}|\nabla v|^2\,dx+ \huno(K\cap \overline A)\,,
$$ 
and the minimality of $v$ yields (\ref{minnA}).
\end{proof}

\begin{proof}[Proof of Theorem \ref{Griffith}.]
Let $t$ be an arbitrary point in $(t_0,t_1)$ and let
$B_i$,  $i=1,\ldots,m$, be a family of open 
balls centred at the points $\gamma_i(\sigma_i(t))$.
If the radii are sufficiently small, we have $\overline 
B_i\subset \Om$ and $\overline B_i\cap K(t_0)=
\overline B_i\cap \overline B_j=\overline B_i\cap 
{\Gamma}_j=\emptyset$ for $j\neq i$. Moreover we may assume 
that $B_i\cap {\Gamma}_i=
\{\gamma_i(\sigma):\tau_i^0<\sigma<\tau_i^1\}$,
for suitable constants $\tau_i^0,\,\tau_i^1$ with 
$\sigma_i^0<\tau_i^0<\sigma_i(t)<\tau_i^1<\sigma_i^1$, and that
the arcs ${\Gamma}_i$ 
intersect $\partial B_i$ only at the points $\gamma_i(\tau_i^0)$ and 
$\gamma_i(\tau_i^1)$, with a transversal intersection. 
All these properties, together with (\ref{structure}), imply that
\begin{equation}\label{Ga}
\overline B_i\cap K(s)=\overline B_i\cap {\Gamma}_i(\sigma_i(s))=
 \{\gamma_i(\sigma):\tau_i^0\leq\sigma\le\sigma_i(s)\}
\quad\hbox{if}\quad \tau_i^0<\sigma_i(s)<\tau_i^1\,.
\end{equation}
In particular this happens for $s=t$, and for $s$ close to $t$ if $\sigma_i$
is continuous at~$t$.

By conditions (a) and (c) of Theorem~\ref{kt} 
and by Lemma~\ref{cond3loc} for every $i$ we 
have that
$$
\E(u(t),K(t),B_i)\le \E(u(t),K,B_i)
\quad\forall \, K\in{\mathcal K}(\overline B_i),\,\ 
K\supset K(t)\cap\overline B_i\,.
$$
By (\ref{Ga}) this implies, taking $K:={\Gamma}_i(\sigma)\cap \overline B_i=
\{\gamma_i(\tau):\tau_i^0\leq\tau\le\sigma\}$,
$$
\E(u(t),{\Gamma}_i(\sigma_i(t)),B_i)\le 
\E(u(t),{\Gamma}_i(\sigma),B_i)\quad \hbox{ for } \sigma_i(t)\le 
\sigma\le \tau_i^1\,,
$$
which yields
\begin{equation}\label{ge0}
\frac{d}{d\sigma} \E(u(t),{\Gamma}_i(\sigma),B_i)
\Big|\lower1.5ex\hbox{$\scriptstyle \sigma=\sigma_i(t)$}
\ge 0\,.
\end{equation}
Inequality (\ref{sif>}) follows now from Lemma~\ref{Grisvard2} applied 
with $g:=u(t)$.

By condition (e) of Theorem \ref{kt} for 
a.e.\ in $t\in (t_0,t_1)$ we have
$\frac{d}{ds}\E(g(t),K(s))|_{s=t}=0$. Moreover, for a.e.\ in $t\in (t_0,t_1)$
the derivative $\dot\sigma_i(t)$ exists 
for $i=1,\ldots,m$. 
Let us fix $t\in (t_0,t_1)$ which satisfies all these properties. 

By (\ref{Ga}) for $s$ close to $t$ we have
\begin{equation}\label{subadd}
\E(g(t),K(s))\le \sum_{i=1}^{m} 
\E(u(t),{\Gamma}_i(\sigma_i(s)),B_i) + \E(u(t),K,A)\,,
\end{equation}
where $K:=K(t_0)\cup \bigcup_i{\Gamma}_i(\tau_i^0)$ and 
$A:=\Om\setmeno \bigcup_i \overline B_i$. Note that the equality holds
in (\ref{subadd}) for $s=t$. As the functions $s\mapsto \E(g(t),K(s))$ 
and $s\mapsto \E(u(t),{\Gamma}_i(\sigma_i(s)),B_i)$ are differentiable 
at $s=t$ (by Theorem~\ref{Grisvard2} and by the existence of 
$\dot\sigma_i(t)$), we 
conclude that
\begin{eqnarray*}
&\displaystyle 
0=\frac{d}{ds} \E(g(t),K(s))
\Big|\lower1.5ex\hbox{$\scriptstyle s=t$}
=\sum_{i=1}^{m} 
\frac{d}{ds} \E(u(t),{\Gamma}_i(\sigma_i(s)),B_i)
\Big|\lower1.5ex\hbox{$\scriptstyle s=t$}=\\
&\displaystyle 
=\sum_{i=1}^{m} \frac{d}{d\sigma} \E(u(t),{\Gamma}_i(\sigma),B_i)
\Big|\lower1.5ex\hbox{$\scriptstyle \sigma=\sigma_i(t)$}
\, \dot\sigma_i(t)=
\sum_{i=1}^{m} \big\{1- \kappa_i(u(t),\sigma_i(t))^2\big\}
\,\dot\sigma_i(t)\,.
\end{eqnarray*}
By (\ref{sigmadot}) and (\ref{sif>}) we have
 $\big\{1- \kappa_i(u(t),\sigma_i(t))^2\big\}
\,\dot\sigma_i(t)\geq0$ for $i=1,\ldots, m$, so that the previous 
equalities yield (\ref{sif=}).
\end{proof}

\end{section}



\begin{thebibliography}{99}

\bibitem{Att}Attouch H.:
Variational Convergence for Functions and Operators.
Pitman, London, 1984.

\bibitem{Att-Pic}Attouch H., Picard C.: Comportement limite de probl\`emes 
de transmission unilateraux a travers des grilles de forme 
quelconque. {\it Rend. Sem. Mat. Politec. Torino\/} {\bf 45} (1987), 
71-85.

\bibitem{BZ} Blake A., Zisserman A.: Visual Reconstruction. MIT 
Press, Cambridge, 1987.

\bibitem{Bou}Bourdin B.: Image segmentation with a finite element 
method. {\it RAIRO Mod\'el. Math. Anal. Num\'er. \/} {\bf 33}  
(1999), 229-244. 


\bibitem{Bou-Cha}Bourdin B., Chambolle A.:
Implementation of an adaptive finite-element 
approximation of the Mumford-Shah functional. {\it Num. Math.\/}
{\bf 85} (2000), 609-646.

\bibitem{BFM}Bourdin B., Francfort G.A., Marigo J.-J.: Numerical 
experiments in revisited brittle fracture. {\it J. Mech. Phys. Solids\/} 
{\bf 48} (2000), 797-826. 

\bibitem{Bre}Brezis H.: Op\'erateurs maximaux monotones et semi-groupes 
de contractions dans les espaces de Hilbert. North-Holland, Amsterdam, 
1973. 


\bibitem{BucVar1}Bucur D., Varchon N.: Boundary variation for the 
Neumann problem. Preprint Univ. Franche-Comt\'e, 
1999. 

\bibitem{BucVar}Bucur D., Varchon N.: A duality approach for the 
boundary variation of Neumann problems. Preprint Univ. Franche-Comt\'e, 
2000. 

\bibitem{Buc}Bucur D., Zolesio J.-P.: $N$-dimensional shape optimization 
under capacitary constraint. {\it J. Differential Equations\/} {\bf 
123} (1995), 504-522. 

\bibitem{Ch} Chambolle A.: Finite-differences discretizations of the 
Mumford-Shah functional. 
{\it RAIRO Mod\'el. Math. Anal. Num\'er. \/} {\bf 33}  (1999), 
261-288. 
 
\bibitem{Ch-D}Chambolle A., Doveri F.: Continuity of Neumann linear 
elliptic problems on varying two-dimensional bounded open sets. {\it Comm. 
Partial Differential Equations\/} {\bf 22} (1997), 811-840. 

\bibitem{ChD}Chambolle A., Doveri F.: Minimizing movements for the 
Mumford and Shah energy. {\it Discrete Continuous Dynamical Systems\/} 
{\bf 3} (1997), 153-174. 

\bibitem{Cor}Cortesani G.: Asymptotic behaviour of a sequence of 
Neumann problems. {\it Comm. Partial Differential Equations\/} {\bf 22} 
(1997), 1691-1729. 




\bibitem{Dam}Damlamian A.: Le probl\`eme de la passoire de Neumann.
{\it Rend. Sem. Mat. Univ. Politec. Torino\/} {\bf 43} (1985), 427-450.

\bibitem{DenLio}Deny J., Lions J.-L.: Les espaces du type de Beppo Levi. 
{\it Ann. Inst. Fourier (Grenoble)\/} {\bf 5} (1953), 305-370. 

\bibitem{ES}Erdogan F., Sih G.C.: On the crack extension in plates 
under plane loading and transverse shear. {\it J. Basic 
Engineering\/} (1963) 519-527.
 
\bibitem{EG}Evans L.C., Gariepy R.F.: Measure Theory and Fine Properties
of Functions. CRC Press, Boca Raton, 1992.



\bibitem{FraMar3}Francfort G.A., Marigo J.-J.: Revisiting brittle 
fracture as an energy minimization problem. {\it J. Mech. Phys. 
Solids\/} {\bf 46} (1998), 1319-1342. 


\bibitem{Gri}Griffith A.: The phenomena of rupture and flow in solids.
{\it Philos. Trans. Roy. Soc. London Ser. A\/} {\bf 221} (1920), 163-198.

\bibitem{Gri1} Grisvard P.: Elliptic Problems in Nonsmooth Domains. 
Pitman, Boston, 1985.

\bibitem{Gri2} Grisvard P.: Singularities in Boundary Value 
Problems. Masson, Paris, 1992.

\bibitem{HKM}Heinonen J., Kilpel\"ainen T., Martio O.: Nonlinear 
Potential Theory of Degenerate Elliptic Equations. Clarendon Press, 
Oxford, 1993. 

\bibitem{Khr}Khruslov E.Ya.: The asymptotic behavior of solutions of the 
second boundary value problem under fragmentation of the boundary of the 
domain. {\it Math. USSR-Sb.\/} {\bf 35} (1979) 266-282.



\bibitem{Ma}Maz'ya V.G.: Sobolev Spaces. Springer-Verlag, Berlin, 
1985. 

\bibitem{MS}Morel J.-M., Solimini S.: Variational Methods in Image 
Segmentation. Birkh\"auser, Boston, 1995.

\bibitem{Mos}Mosco U.: Convergence of convex sets and of solutions 
of variational inequalities. {\it Adv. in Math.\/} {\bf 3} (1969), 
510-585.

\bibitem{MSh} Mumford D., Shah A.: Optimal approximation by 
piecewise smooth functions and associated variational problems. 
{\it Comm. Pure Appl. Math.\/} {\bf 42} (1989), 577-685.

\bibitem{Mur}Murat F.: The Neumann sieve. {\it Non Linear Variational 
Problems (Isola d'Elba, 1983)\/}, 27-32, {\it Res. Notes in Math. 127, 
Pitman, London\/}, 1985.

\bibitem{Rich}Richardson T.J.: Scale independent piecewise 
smooth segmentation of images via variational methods. 
Ph.D. Thesis, Department of Electrical Engineering 
and Computer Science, MIT, Cambridge (MA, USA), 1990.

\bibitem{RM}Richardson T.J., Mitter S.K.: A variational 
formulation-based edge focussing algorithm. {\it Sadhana} {\bf 22} 
(1997), 553-574.

\bibitem{Rog} Rogers C.A.: Hausdorff Measures. Cambridge University 
Press, Cambridge, 1970.

\bibitem{SL} Sih G.C., Liebowitz H.: Mathematical theories of brittle 
fracture. {\it Fracture: An Advanced Treatise, vol. II, Mathematical 
Fundamentals, ed. Liebowitz H.\/}, 67-190, {\it Academic Press, New 
York\/}, 1968.

\bibitem{SM} Sih G.C., Macdonald B.: Fracture mechanics applied to 
engineering problems -- strain energy density fracture criterion. {\it 
Engineering Fracture Mechanics\/} {\bf 6} (1974) 361-386. 

\bibitem{Sve}\v Sverak V.: On optimal shape  design. 
{\it J. Math. Pures  Appl.\/} {\bf 72} (1993), 537-551. 

\bibitem{Zie}Ziemer W.P.: Weakly Differentiable Functions. 
Springer-Verlag, Berlin, 1989. 

\end{thebibliography}
\end{document}